\documentclass[12pt,a4paper]{article}
\usepackage{preamble}
\usepackage{vmargin}
\setmargrb{1in}{0.3in}{1in}{0.7in} 
\newcommand{\usualHilb}{{\ccH}ilb}
\newcommand{\minD}[1][D]{\operatorname{min}_{#1}}
\DeclareMathOperator{\sat}{sat}

\DeclareMathOperator{\QDSG}{QDSG}

\definecolor{zielony}{rgb}{0.2, 0.5, 0}
\definecolor{czerwony}{rgb}{0.9, 0.2, 0.1}
\definecolor{brazowy}{rgb}{0.5, 0.1, 0.0}
\definecolor{niebieski}{rgb}{0.3, 0.1, 0.9}
\definecolor{skin}{HTML}{FFECC9}
\definecolor{pumpkin}{HTML}{FEDFA9}
\definecolor{piggy}{HTML}{FFB99D}
\definecolor{fiolet}{HTML}{CD8F9C}
\definecolor{granat}{HTML}{677081}
\definecolor{ciemnyblekit}{HTML}{91A1B8}

\title{Cactus varieties of sufficiently ample embeddings of projective schemes
have determinantal equations\thanks{Dedicated in memory of Gianfranco Casnati}}
\author{\nisiabu \and \JaBu \and {\L}ucja Farnik}

 \date{30 November 2024}
\newenvironment{red}{\color{red}}{}
\newcommand{\bred}{\begin{red}}
\newcommand{\ered}{\end{red}}

\newenvironment{blue}{\color{blue}}{}
\newcommand{\bblue}{\begin{blue}}
\newcommand{\eblue}{\end{blue}}

\newenvironment{green}{\color{green}}{}
\newcommand{\bgreen}{\begin{green}}
\newcommand{\egreen}{\end{green}}

\usepackage[textsize=tiny]{todonotes}

\newcommand{\Apolar}[1][p]{\operatorname{Apolar}(#1)}

\DeclareMathOperator{\cone}{cone}

\DeclareMathOperator{\Amp}{Amp}
\DeclareMathOperator{\interior}{int}

\definecolor{DziubkaDyniowy}{RGB}{224,104,19}
\definecolor{geometria}{RGB}{170,67,6}
\definecolor{DziubkaBrazA}{RGB}{140,64,3}
\definecolor{algebra}{RGB}{75,135,10} 
\begin{document}
\maketitle
\begin{abstract}
For a fixed projective scheme $X$, a property $P$ of line bundles is satisfied by sufficiently ample line bundles
if there exists a line bundle $L_0$ on $X$ such that $P(L)$ holds for any $L$ with $L-L_0$ ample.
As an example, sufficiently ample line bundles are very ample, moreover, for a normal variety $X$,
the embedding corresponding to sufficiently ample line bundle is projectively normal.
The grandfather of such properties and a basic ingredient used to study this concept is Fujita vanishing theorem,
which is a strengthening of Serre vanishing to sufficiently ample line bundles.
The $r$-th cactus variety of $X$ is an analogue of secant variety and it is defined using
linear spans of finite schemes of degree $r$. In this article we show that cactus varieties of sufficiently ample embeddings
of $X$ are set-theoretically defined by minors of matrices with linear entries.
The topic is closely related to conjectures of Eisenbud-Koh-Stillman, which was proved by Ginensky
in the case  $X$  a smooth curve~\cite{ginensky_Eisenbud_for_curves}.
On the other hand Sidman-Smith~\cite{sidman_smith_linear_det_eq_for_all_proj_schemes}
proved that the ideal of sufficiently ample embedding of  any projective scheme $X$ is
generated by $2 \times 2$ minors of a matrix with linear entries.
\end{abstract}

\medskip
{\footnotesize
\noindent\textbf{addresses:} \\
J.~Buczy\'nski, \eemail{jabu@mimuw.edu.pl}, 
Institute of Mathematics of the Polish Academy of Sciences, ul. \'Snia\-deckich 8, 00-656 Warsaw, Poland\\
W.~Buczy\'nska, \eemail{wkrych@mimuw.edu.pl}, 
Faculty of Mathematics, Computer Science and Mechanics, University of Warsaw,
ul. Banacha 2, 02-097 Warsaw, Poland\\
{\L}.~Farnik, \eemail{lucja.farnik@uken.krakow.pl},
Department of Mathematics, University of the National Education Commission, ul. Podchorążych 2, 30-084 Kraków, Poland\\

\noindent\textbf{keywords}:
cactus varieties, sufficiently ample line bundles, determinantal equations, finite schemes

}

\tableofcontents

\section{Introduction}
Throughout most of this article we work over a fixed algebraically closed field $\kk$ of any characteristics.
Only in Appendix~\ref{sec_appendix_nonclosed_fields} we use more general base fields.

\subsection{Secant and cactus varieties}
\label{sec_intro_secant_and_cactus}

Suppose $X$ is a projective variety over $\kk$ and  $r$ is a fixed  positive integer.
For a closed embedding $X\subset \PP^N$,
the $r$-th secant variety $\sigma_r(X)$ of $X$ (strictly speaking, of $X\subset \PP^N$),   is defined as
$\sigma_r(X) = \overline{\bigcup \set{ \linspan{\fromto{p_1}{p_r}}
\mid p_i\in X}} \subset \PP^N$.
Here, $\setfromto{p_1}{p_r} \subset X$ is a smooth finite subscheme,
and $\linspan{\dots}$ denotes the (projective) linear span of a subscheme in $\PP^N$.
Secant varieties have practical applications both in algebraic geometry and also outside this area of mathematics,
see for instance \cite[Chapt.~1]{landsberg_tensorbook},
\cite[Chapt.~1]{landsberg_geometry_and_complexity}.
It is a long standing problem to find the defining equations for secant varieties. The first and most natural equations
that vanish on the $r^{\text{th}}$ secant variety of the Veronese variety are $(r+1) \times (r+1)$ minors
of the catalecticant matrices.
More generally, elementary linear algebra and semicontinuity of the matrix rank show that if
$X\subset \set{\rk M \le k}$ for some matrix $M$ with entries that are linear coordinates on $\PP^N$,
then $\sigma_r(X)\subset \set{\rk M \le k\cdot r}$.
By means of  this observation, more families of equations were found, see a summary in
in~\cite{landsberg_ottaviani_VB_method_equ_for_secants}.
A natural method introduced there produces a matrix with linear entries from vector bundles on $X$.
Unfortunately, the extrinsic geometry of secant varieties in $\PP^N$ tends to be very difficult and behaves wildly,
except in a few tame cases.

In \cite{nisiabu_jabu_cactus}, \cite{galazka_vb_cactus}, and
\cite[Thm~1.18]{galazka_phd} the wild beast is given a name by exhibiting a potentially larger,
typically reducible, variety, called \emph{cactus variety}  $\cactus{r}{X}$, defined below in \eqref{equ_define_cactus_variety}.
Again, strictly speaking, the cactus variety depends not only on $X$, but also on its embedding
$\iota\colon X \to \PP^N$.
But as in the case of the secant varieties, unless potentially confusing, it is customary to omit $\iota$ from the notation.
Moreover, in the definition it is enough to assume that $X$ is a projective scheme, not necessarily a variety.
\begin{equation}
\label{equ_define_cactus_variety}
\cactus{r}{X}:= \overline{\bigcup \set{\langle R \rangle  \mid R \subset X,
\text{ $R$ is a subscheme of degree at most $r$}  }} \subset \PP^N.
\end{equation}
The point is that (if $X$ is a projective variety) on one hand $\sigma_r(X) \subset \cactus{r}{X}$ and on
the other hand many natural equations vanishing on $\sigma_r(X)$ vanish also on the larger set
$\cactus{r}{X}$.
Thus the difference $\cactus{r}{X} \setminus \sigma_r(X)$ is a natural, geometric obstruction to the ``tameness''
of secant varieties. The goal of this article is to show that under certain assumptions this difference is the only obstruction.
We show that for sufficiently ample embedding $X\subset \PP^N$ its cactus variety $\cactus{r}{X}$ is set-theoretically
defined by minors of a matrix, whose entries are linear forms on $\PP^N$. These minors arise naturally from
sections of certain line bundles.

\subsection{Cactus variety of a sufficiently ample embedding}

For the purpose of this introduction we define the properties satisfied by ``sufficiently ample line bundles'' on $X$
in the following way, and we further discuss this notion in Section~\ref{sec_suff_ample_line_bndles}.

\begin{defin}
\label{def_sufficiently_ample_l_b}
Fix a projective scheme $X$ over $\kk$.
Let $P$ be a property of line bundles on $X$.
We say that $P$ is satisfied for sufficiently ample line bundles
if there exists a line bundle $L_0$ such that $P$ holds for all $L_0 \otimes L$ with $L$ ample.
\end{defin}
In particular, ``very ampleness'' is satisfied by sufficiently ample line bundles,
which in less strict, but more natural language we phrase as
``a sufficiently ample line bundle on $X$ is very ample''.
See Section~\ref{sec_suff_ample_line_bndles} for discussion, references and a concise proof of this fact.

Our main result, in a simplified form, is the following.
\begin{thm}\label{thm_main_intro_simpler}
Fix an integer $r'$, and suppose $X$ is a projective scheme.
For a sufficiently ample line bundle $L$
consider embedding of $X$ given by the complete linear system $H^0(X,L)$:
\[
\iota \colon X \hookrightarrow \PP(H^0(X,L)^*) \simeq \PP^N.
\]
Then there exists a matrix $M$ whose entries are
are linear coordinates on $\PP^N$, equivalently they are sections of $L$
-- elements of $H^0(X, L)$, such that
\[
\cactus{r}{X} = \reduced{\set{\rk M \le r}} \quad  \text{ for all } r\le r',
\]
where $\reduced{\set{\rk M \le r}}$ is the (reduced) zero set of the homogeneous ideal generated by
$(r+1)\times (r+1)$ minors of $M$.
\end{thm}

In particular, the theorem gives a partial answer to Conjecture of
Sidman and Smith, see~\cite[Conj.~1.2]{sidman_smith_linear_det_eq_for_all_proj_schemes}.
We discuss this and similar conjectures together with their relations to our results
in \S\ref{sec_conjectures_secants_to_suff_ample}.

More precisely, the matrix $M$ in Theorem~\ref{thm_main_intro_simpler}
generalises the catalecticant matrix in the case of $X=\PP^n$, compare to \cite[Thm~1.5]{nisiabu_jabu_cactus}.
Given a splitting $L= A\otimes B$ as a product of  line bundles $A$ and $B$ on $X$,
there is a corresponding map --- multiplication of sections: $H^0(X, A)\otimes H^0(X, B) \to H^0(X, L)$.
For any choice of bases $\fromto{s_1}{s_k}$ of $H^0(X, A)$ and $\fromto{t_1}{t_l}$ of $H^0(X, B)$
we obtain a $k\times l$ matrix $M$, whose $(i,j)$-th entry is $s_i t_j\in H^0(X, L)$.
We claim that  when both $A$ and $B$ are sufficiently ample and satisfy additional technical conditions,
in some sense, $A$ and $B$ as divisors are close to a ``half'' of $L$,
the matrix $M$ as above has the desired properties of Theorem~\ref{thm_main_intro_simpler}.

\subsection{Conjectures on secant varieties to sufficiently ample embeddings}
\label{sec_conjectures_secants_to_suff_ample}

In a series of results by Mumford  \cite[p.~32, Thm~1]{mumford_varieties_defined_by_quadrics},
Griffiths \cite[Thm~p.~271]{griffiths_infsmal_var_Hodge_structure_3det_varieties},
and Sidman and Smith \cite[Thm~1.1]{sidman_smith_linear_det_eq_for_all_proj_schemes}
it is shown that for any projective scheme $X$ the ideal of a sufficiently ample embedding of $X$
is generated by $2\times 2$ minors of a matrix with linear entries, specifically, one constructed as above.
In particular,
\cite[Thm~1.1]{sidman_smith_linear_det_eq_for_all_proj_schemes}
prove a version of Theorem~\ref{thm_main_intro_simpler} for $r'=1$.
Moreover they phrase the following conjecture \cite[Conj~1.2]{sidman_smith_linear_det_eq_for_all_proj_schemes}.
A version of this conjecture for curves is stated also by
Eisenbud, Koh, and Stillman  \cite[p.~518, Equation~(*)]{eisenbud_koh_stillman_curves_of_high_degree}.

\begin{conjecture} \label{conj_Sidman_Smith_and_E_K_S}
Fix an integer $r$ and  a projective variety  $X$.
Suppose $L$ is a sufficiently ample line bundle on $X$,
with $\iota \colon X\to \PP(H^0(X, L))$ the embedding by the complete linear system.
Then the ideal of $\sigma_r(X)$ is generated by $(r+1)\times (r+1)$ minors of the matrix with linear entries.
\end{conjecture}

Conjecture \ref{conj_Sidman_Smith_and_E_K_S} was proved set-theoretically  in the case $X$ is a smooth curve
in  \cite{ravi_det_equations_for_sec_vars_of_curves},  and scheme-theoretically for smooth curves
in  \cite{ginensky_Eisenbud_for_curves}.
On the other hand, if $X$ has  sufficiently bad singularities, even if it  is a curve,
then \cite[Thm~1.17]{jabu_ginensky_landsberg_Eisenbuds_conjecture}
provides counterexamples to Conjecture~\ref{conj_Sidman_Smith_and_E_K_S}.
At this point, it should not come as  a surprise to the reader, that the real reason for the failure of this conjecture
comes from our favourite obstruction, that is the difference between cactus and secant varieties:
$\cactus{r}{X} \setminus\sigma_r(X)$.
For this reason, in the paragraph following
the conjecture~\cite[Conj.~1.2]{sidman_smith_linear_det_eq_for_all_proj_schemes}
the authors mention that ``\cite[Thm~1.3]{nisiabu_jabu_cactus} suggests that the secant varieties
in Conjecture~\ref{conj_Sidman_Smith_and_E_K_S} should be replaced by cactus varieties''.
This is exactly what we prove in Theorem~\ref{thm_main_intro_simpler}, at the moment in the set-theoretic version.

\begin{rem}
The article \cite{sidman_smith_linear_det_eq_for_all_proj_schemes}
contains additional thesis on the resulting matrix $M$, both in Theorem~1.1, and in Conjecture~1.2.
Namely, the authors claim that $M$ is $1$-generic, that is the orbit of $M$ with respect to
the action of the product of $\GL(H^0(X,A))\times\GL(H^0(X,B))$ does not contain any matrix with a zero entry.
If $X$ is a variety (reduced and irreducible), then this condition is automatic from the construction of the matrix
as one arising from the multiplication of sections, and from \cite[Prop.~6.10]{eisenbud_geometry_of_syzygies}.
Otherwise, if $X$ is not irreducible or not reduced,
then $1$-genericity claim in~\cite[Thm~1.1]{sidman_smith_linear_det_eq_for_all_proj_schemes} fails.
For example, if $X=\Spec \kk[t]/(t^2-t)$ or $X=\Spec \kk[t]/(t^2)$, then there is no $1$-generic matrix with linear entries,
such that the ideal of $X$ is generated by its $2\times 2$ minors in any projective embedding.
Therefore, in our considerations we ignore the issue of $1$-genericity of $M$.
\end{rem}

There exist some cases when the cactus and secant varieties are equal. In those situations  our results
provide set-theoretical equations of secant varieties.
\begin{cor}
Suppose $X$ is a projective variety of dimension $n$ and $r'$ is an integer.
If one of the following conditions hold:
\begin{itemize}
\item $X$ is smooth and $n\leqslant 3$, or
\item $X$ is smooth and $\cchar \kk \neq 2, 3$ and $r'\leqslant 13$, or
\item $X$ is smooth and $\cchar \kk \neq 2, 3$, $r'= 14$, and $n\leqslant 5$, or
\item $n=1$ and $X$ is a curve with only planar singularities,
\end{itemize}
then for a sufficiently ample embedding of $X$, in the notation of Theorem~\ref{thm_main_intro_simpler},
the secant variety $\sigma_r(X)$ is set theoretically defined by $(r+1)\times (r+1)$ minors of $M$ for all $r\leqslant r'$.
\end{cor}
This corollary is an immediate consequence of Theorem~\ref{thm_main_intro_simpler} combined
with~\cite[Cor.~6.20, \S7.2]{jabu_jelisiejew_finite_schemes_and_secants}. In particular, the explicit conditions
on dimension and $r$ are copied from Tables~1--4 in~\cite{jabu_jelisiejew_finite_schemes_and_secants}.
Note that there is some space for improving the number of cases, where secant variety and cactus variety
coincide in characteristic $2$ and $3$: one can also show that in any characteristic
$\sigma_r(X) = \cactus{r}{X}$ for $r\leqslant 5$ and $X$ non-singular.

Moreover, fine properties of secant varieties to sufficiently ample embeddings has attracted recently
a significant attention, so far mainly for smooth projective curves, as this is the most tractable case, and also
the first testing ground:
\cite{ein_niu_park_sings_and_syz_of_sec_vars_of_curves},
\cite{choe_kwak_park_syzygies_sec_var_of_sm_proj_curves}.
An ongoing independent work  in preparation \cite{choi_lacini_park_sheridan_sings_and_syz_of_secant_vars}
also investigates cases of surfaces and more general varieties.

\subsection{Overview}
In Section~\ref{sec_cones} we study combinatorial properties of not necessarily polyhedral cones
in real affine spaces with their relations to lattice points.
This serves as a preparation for a discussion of $\Nef(X)$-cone and its relation to Picard group.
In Section~\ref{sec_suff_ample_line_bndles},
we  present basic facts about sufficiently ample line bundles.
In Section~\ref{sec_multigraded_rings} we study
multigraded rings coming from sections of line bundles in a fairly restricted setting, commencing with
a single-graded ring of sections of (powers of) a line bundle and then generalising to a double-graded one.
In particular, in \S\ref{sec_apolarity} we discuss a version of apolarity in this setting.
In Section~\ref{sec_cacti} we prove our main result. The proof is divided into two parts: algebraic, where
we exploit the apolarity, and combinatorial, where we show that some convex body must have a lattice point,
which then leads to an appropriate splitting of a sufficiently ample line bundle $L=A\otimes B$.
In Section~\ref{sec_future} we briefly mention the further or ongoing projects related to this topic.
Finally, in  Appendix~\ref{sec_appendix_nonclosed_fields}
we provide references and sketch the argument that generalises our result to arbitrary base fields.

\subsection*{Acknowledgements}
We are grateful to Piotr Achinger, Joachim Jelisiejew, Wojciech Kucharz, Alex K{\"u}ronya,
Gregory G.~Smith, Frank Sottile, Tomasz Szemberg, and Jarosław Wiśniewski
for interesting discussions and helpful comments.
The authors are supported by three research grants awarded by National Science Center, Poland:
W.~Buczy{\'n}ska and J.~Bu\-czy{\'n}\-ski are supported by
the project ``Complex contact manifolds and geometry of secants'', 2017/26/E/ST1/00231.
{\L}.~Farnik is supported by the project ``Positivity of line bundles on algebraic surfaces'', 2018/28/C/ST1/00339.
In addition, J.~Buczy{\'n}ski is partially supported by the project ``Advanced problems in contact manifolds,
secant varieties, and their generalisations (\mbox{APRICOTS+})'',  2023/51/B/ST1/02799.
Moreover, part of the research towards the results of this article was done during the scientific semester
Algebraic Geometry with Applications to Tensors and Secants in Warsaw.
The authors are grateful to many fruitful discussions with the participants, and for the partial support
by the Thematic Research Programme ``Tensors: geometry, complexity and quantum entanglement'', University of Warsaw,
Excellence Initiative---Research University and the Simons Foundation Award No.~663281 granted to
the Institute of Mathematics of the Polish Academy of Sciences for the years 2021-2023.

\section{Cones, lattices and numerical equivalence of divisors}\label{sec_cones}

We gather some general results of convex geometry flavour that we  need to prove
some statements about nef and ample cones of schemes.

\subsection{Lattice cones}
Consider a lattice $N= \ZZ^{\rho}$
and its underlying $\RR$-vector space
$N_{\RR} = N\otimes_{\ZZ} \RR \simeq \RR^{\rho}$.
On $N_{\RR}$ we  consider Euclidean topology: for a subset $\sigma\subset N_{\RR}$ by $\overline \sigma$ we denote
its closure and by $\interior(\sigma)$ or $\interior_N(\sigma)$ the  interior of  $\sigma$ inside $N_{\RR}$.
We always have the natural inclusion $N \subset N_{\RR}$ arising from $\ZZ\subset \RR$.

A subset $\sigma \subset N_{\RR}$ is called
a \emph{cone} if it is convex and invariant under the positive rescalings: $\RR_{> 0}\cdot \sigma = \sigma$.
For any subset $S\subset N$, we let $\cone(S)$ to be the smallest cone in $N_{\RR}$ containing $S$,
or equivalently:
\[
\cone(S)=\set{\textstyle{\sum_{i=1}^k} c_i s_i \mid k\in \ZZ, c_i\in \RR_{>0}, s_i \in S}.
\]
Every cone has its \emph{dimension} defined as the $\RR$-dimension of its $\RR$-linear span.
For two subsets $S_1, S_2 \subset N_{\RR}$ their \emph{Minkowski sum} $S_1+S_2$ is the set
$\set{s_1 + s_2\mid s_1 \in S_1, s_2 \in S_2}$.
For $D\in N_{\RR}$ and a cone $\sigma\subset N_{\RR}$, the sum $\set{D}+\sigma$
(denoted simply $D+\sigma$) is a \emph{translated cone}.

Our main interest is in cones $\sigma\subset N_{\RR}$ of \emph{full dimension},
which are characterised by one of the following equivalent conditions:
\begin{itemize}
\item the $\RR$-linear span of $\sigma$ is $N_{\RR}$,
\item $\interior(\sigma)\neq \emptyset$,
\item if $V\subset N_{\RR}$ is the linear span of $\sigma$, then for any $L\in \interior_V(\sigma)$ we have
$\cone(-L) + \sigma = N_{\RR}$.
\end{itemize}

\begin{lemma}
\label{lem_cone_and interior_lattice_points}
Suppose $\sigma\subset N_{\RR}$
is a cone of full dimension, and $P\subset N$ is any subset.
Then the following two conditions are equivalent:
\begin{itemize}
\item $D+\sigma\cap N  \subset P$ for some $D\in N$,
\item $D'+\interior(\sigma)\cap N \subset P$ for some $D'\in N$.
\end{itemize}
\end{lemma}
\begin{prf}
Assuming the first item, let $D'=D$ and the second property is automatically satisfied.
For any $L\in \interior(\sigma)\cap N$,
we have $L+\sigma \subset \interior(\sigma)$.
Thus assuming the second item, take $D=D'+L$,
and the first property holds.
\end{prf}

\begin{lemma}
\label{lem_two_cones_in_lattice_and_sublattice}
Suppose  we have two lattices  $N'\subset N$ and cones
$\sigma\subset N_{\RR}$ and $\sigma'\subset N'_{\RR}$ that are  of full dimensions
in the respective  lattice.
Assume $\sigma'\subset \sigma$ and  $\interior_{N'}(\sigma') \subset \interior_{N}(\sigma)$.
Then for any $D\in N$ we have
$N'\cap (D+\sigma\cap N) \ne \emptyset$
and there exists $D'\in N'$ such that
$D'+\sigma' \subset D+ \sigma$.
\end{lemma}

\begin{prf}
The interior of $\sigma'$ in $N'_{\RR}$   must contain a lattice point $L'\in N'$,
which is also an interior point of $\sigma$.
Thus $N_{\RR}=\cone(-L')+\sigma \supset -D+N'$,
and hence for any $E'\in N'$ for some integer $k\ge 0$
we must have $E'+ k \cdot L' \in (D+ \sigma) \cap N'$.
Therefore  $N'\cap (D+(\sigma\cap N)) \ne \emptyset$ as claimed.

To show the second claim pick any $D'\in N'\cap (D+(\sigma\cap N))$.
Then $D'+\sigma' \subset D'+\sigma \subset D+\sigma$ as claimed.
\end{prf}

When $\sigma$ is a closed and pointed cone of full dimension,
we define the \emph{dual cone} of $\sigma$
\[
\sigma^{\vee}
:=
\set{\phi \in N_{\RR}^* \mid \phi(D)\geqslant 0 \text{ for all } D \in \sigma}.
\]
Since $\sigma$ is a closed pointed cone of full dimension, also $\sigma^{\vee}$ is a closed pointed cone of full dimension.
For any $H\in \interior(\sigma)$
consider the corresponding hyperplane $(H=1)$ in $N_{\RR}^*$:
\[
(H=1) := \set{\phi \in N^1_{\RR}(X)^* \mid \phi(H)=1}.
\]
Then in the Euclidean topology $\sigma^{\vee}\cap (H=1)$ is a compact subset of $N_{\RR}^*\simeq \RR^{\rho}$.
Moreover, $\sigma^{\vee}\cap (H=1)$  generates the cone:
$\sigma^{\vee} = \RR_{\geqslant 0 }\cdot \left(\sigma^{\vee} \cap (H=1)\right)$.
Thus by double duality
\begin{subequations}
\begin{align}
\sigma &=
\set{D \in N_{\RR} \mid \phi(D)\geqslant 0 \text{ for all } \phi \in \sigma^{\vee} \cap (H=1)}, \text{ and}
\label{equ_closed_cone_by_its_inequalities}\\
\interior(\sigma) &=
\set{D \in N_{\RR} \mid \phi(D)> 0 \text{ for all } \phi \in \sigma^{\vee} \cap (H=1)}.
\label{equ_open_cone_by_its_inequalities}
\end{align}
\end{subequations}

\begin{lemma}
\label{lem_minD_and_maxD}
Suppose $\sigma$ is a closed and pointed cone of full dimension and $H$ is as above
For any $D\in N_{\RR}$
there exists real number
\[
\minD = \min \set{\phi(D) \mid \phi\in \sigma^{\vee} \cap (H=1)}.
\]
In addition:
\begin{itemize}
\item $D\in \sigma$ if and only if $\minD\geqslant 0$, and
\item  $D\in \interior(\sigma)$ if and only if $\minD>0$.
\end{itemize}
\end{lemma}
\begin{prf}
Since $\sigma^{\vee} \cap (H=1)$ is compact, the continuous function $\phi \mapsto \phi(D)$ attains its minimum $\minD$.
The additional claims follow from \eqref{equ_closed_cone_by_its_inequalities}
and \eqref{equ_open_cone_by_its_inequalities}.
\end{prf}

Now fix $\sigma$ and $H$ as above and pick any basis $\fromto{e_1}{e_{\rho}}$ of $N$.
We will consider a closed cube $C\subset N_{\RR}$
of size $2$, 
whose vertices are $\pm e_1 \pm \dotsb \pm e_{\rho}$.

\begin{lemma}
\label{lem_condition_on_lattice_point_in_cones}
Suppose $D'$ and $D''$ are two vectors in $N_{\RR}$.
If for any choice of a vertex $v=\pm e_1 \pm \dotsb \pm e_{\rho}$ of the cube $C$
we have
$D' -D'' + v \in \sigma$,
then there exists a lattice point in
\[
(D' - \sigma ) \cap (D''+ \sigma) \cap N.
\]
\end{lemma}
\begin{prf}
There is always a lattice point in every translation of $\frac{1}{2} C$ (a cube of size $1$).
We claim that
\begin{equation}
\label{equ_condition_on_cube}
\tfrac{1}{2}(D'+D'') + \tfrac{1}{2} C \subset (D' - \sigma ) \cap (D''+ \sigma),
\end{equation}
thus the intersection contains a lattice point.
The set $(D' - \sigma ) \cap (D''+ \sigma)$ is convex so it is enough to check the containment on vertices of $C$ only,
which we do below. Let $v=\pm e_1 \pm \dotsb \pm e_{\rho}$ be any vertex of $C$.
By our assumption $D' -D'' + v \in \sigma$.
Since $\sigma$ is a cone (invariant under positive rescaling), this is equivalent to
\begin{align*}
\tfrac{1}{2}(D' -D'') + \tfrac{1}{2}v &\in \sigma.\\
\intertext{
Subtracting $D'$ or adding $D''$ to both sides the above is equivalent to either of the following:
}
\tfrac{1}{2}(D' + D'') + \tfrac{1}{2}v &\in D''+\sigma, \text{ or}\\
\tfrac{1}{2}(D' + D'') - \tfrac{1}{2}v &\in D' - \sigma.
\end{align*}
The last two conditions imply \eqref{equ_condition_on_cube} and hence conclude the proof.
\end{prf}
\subsection{Ample and nef cones}\label{sec_ample_and_nef_cones}
Suppose $X$ is a projective scheme over $\kk$ and let $\Pic(X)$ denote the Picard group of $X$.
We say that two line bundles $L_1$ and $L_2$ are \emph{numerically equivalent} if for any irreducible and reduced
curve $C\subset X$ the restricted line bundles $L_1|_C$ and $L_2|_C$ have the same degree.
Let $N^1(X): =\Pic(X) / \equiv$, where $\equiv$ is the numerical equivalence.
By \cite[Prop~IV.1.4]{kleiman_numerical_theory_of_ampleness}
see also \cite[Thm~1.1.16]{PAG_I} for more modern treatment but to some extent restricted to $\kk=\CC$,
we have $N^1(X)\simeq \ZZ^{\rho}$ for some integer $\rho$ called the \emph{Picard number} of $X$.

By the Nakai-Moishezon-Kleiman criterion
\cite[Thm~III.1.1(i) and (ii)]{kleiman_numerical_theory_of_ampleness},
see also \cite[Thm~1.2.23 and Cor.~1.2.24]{PAG_I},
if $L_1$ and $L_2$ are numerically equivalent and $L_1$ is ample then $L_2$ is ample.
Therefore it makes sense to define the set of ample classes $\Amp_{\ZZ}(X)\subset N^1(X)$.
The ample cone $\Amp(X)\subset N_{\RR}^1(X) = N^1(X)\otimes_{\ZZ} \RR \simeq \RR^{\rho}$
is then defined as  $\cone(\Amp_{\ZZ}(X))$.
By the Nakai-Moishezon-Kleiman criterion for $\RR$-divisors
\cite[Thm~1.3]{campana_peternell_algebraicity_of_ample_cone} we also have $\Amp(X)$
described by the same intersection theory inequalities as $\Amp_{\ZZ}(X)$,
and thus $\Amp_{\ZZ}(X) = \Amp(X)\cap N^1(X)$.

A line bundle $L$ is \emph{nef} if its restriction $L|_C$ has non-negative degree for every irreducible and reduced
curve $C\subset X$. A numerical class $[L]$ in $N^1(X)$ is \emph{nef} if its representative
$L\in \Pic(X)$ is nef, which does not depend on the choice of the representative.
The set of nef classes is denoted
$\Nef_{\ZZ}(X)$ and the corresponding cone  (called the \emph{nef cone}) in $N_{\RR}^1(X) $ is
$\Nef(X) = \cone\left(\Nef_{\ZZ}(X)\right)$.

Since $X$ is projective, we have
\begin{equation}\label{equ_interior_and_closure_Amp_and_Nef}
\Nef(X) = \overline{\Amp(X)},
\quad\text{and} \quad \interior (\Nef(X)) = \Amp(X).
\end{equation}
In particular, $\Amp(X)$ is open and non-empty by projectivity, thus both cones have full dimension.

Let $\overline{NE(X)} \subset N^1_{\RR}(X)^*$ be the \emph{dual cone} of $\Nef(X)$:
\[
\overline{NE(X)} :=
\set{\phi \in N^1_{\RR}(X)^* \mid \phi(D)\geqslant 0 \text{ for all } D \in \Nef(X)}.
\]
If $X$ is a projective variety, then the dual cone is traditionally called the \emph{Mori cone},
or the closure of the \emph{cone of curves}, and denoted $\overline{NE(X)}$.
We extend this notation to our more general situation of $X$ being the projective scheme.

An immediate consequence of Lemma~\ref{lem_cone_and interior_lattice_points}
and \eqref{equ_interior_and_closure_Amp_and_Nef} is the following proposition.
\begin{prop}
\label{prop_nef_definition_of_sufficiently_ample}
Let $X$ be projective scheme and consider
a property $P$ of line bundles on $X$.
Then the following are equivalent:
\begin{enumerate}
\item $P$ is satisfied for sufficiently ample line bundles in the sense of Definition~\ref{def_sufficiently_ample_l_b},
\item there exists a class $D\in N^1(X)$ such that $P$ holds for any line bundle $L$ whose numerical class is
  in contained in the translated cone $\sigma := D + \Nef(X)$.
\end{enumerate}
\end{prop}
\noprf

Trivial examples of properties that are
satisfied for sufficiently ample line bundles are \emph{ampleness} and \emph{nefness}.
Another example is provided by
Fujita vanishing (Corollary~\ref{cor_fujita_vanishing}).
Moreover, sufficiently ample line bundles are
\emph{very ample} (Proposition~\ref{prop_sufficiently_ample_is_very_ample}).

Note that the starting divisor $D$ is not uniquely determined, and a priori we do not require that $D$ itself is nef or ample.
However, it is always possible to choose $D$ nef or even ample by replacing any $D$ by $D+d\Delta$
for a sufficiently large multiple $d$, where $\Delta$ is any ample divisor.

In this article we do not discuss ``optimal'' or ``best'' choices of $D\in N^1(X)$ such that
a given property $P$ holds for all line bundles whose numerical classes are in $D+\Nef(X)$.
We are content with statements just that such $D$ exists for the properties considered here.
However, the problem of describing or partially describing the set of all the line bundles that satisfy
a given property is also very interesting and should be addressed in future research.
We only discuss a single, elementary example in Example~\ref{ex_hypereliptic_curve}.
We also prove Theorem~\ref{thm_main_with_splitting}, a version of our main result,
which explicitly lists properties required to obtain the description of cactus varieties as zeroes of minors.
This could be a base for future research towards an effective version of Theorem~\ref{thm_main_intro_simpler}.

\begin{lemma}
\label{lem_finitely_many_properties_of_sufficiently_ample_lbs}
Suppose $\fromto{P_1}{P_k}$ are properties of line bundles on a fixed projective scheme $X$ and assume
that each $P_i$ is  satisfied for sufficiently ample line bundles.
Then for sufficiently ample line bundles $L$, all properties $\fromto{P_1}{P_k}$ hold.
\end{lemma}

\begin{prf}
  For each property $P_i$ there exists a divisor class $D_i$ such that $P_i$ is satisfied for all line bundles
  with numerical class in $D_i+\Nef(X)$.
Without loss of generality we may assume that each $D_i$ is nef.
Then for $D= D_1+ \dotsb + D_k$ we have
$D+\Nef(X) \subset D_i + \Nef(X)$ and thus for
line bundles whose classes are in $D+\Nef(X)$ all the properties $P_i$ hold.
\end{prf}

\begin{example}
\label{ex_infinitely_many_properties_of_suff_ample_lbs}
The conclusion of Lemma~\ref{lem_finitely_many_properties_of_sufficiently_ample_lbs} is incorrect if
we insist on infinitely many properties. As an immediate example, for $X=\PP^1$ and each $i \in \ZZ$
let $P_i$ be the property of a line bundle $L$ on $\PP^1$ defined as:
\[
P_i (L) \iff H^1(L\otimes \ccO_{\PP^1}(i)) =0.
\]
Then for each $i$ and a sufficiently ample line bundle $L$ the property $P_i$ holds.
However, there is no line bundle $L$ that has all properties $P_i$ for every $i$.
\end{example}
\subsection{Structure of Picard group of a product}

For two schemes $X_1$ and $X_2$ and line bundles $L_i$ on $X_i$ (for $i=1,2$),
we  consider the product $X_1\times X_2$ with two projections $\pi_i \colon X_1\times X_2 \to X_i$.
Define the product line bundle $L_1\boxtimes L_2$ on $X_1\times X_2$ as
\[
L_1\boxtimes L_2:= \pi_1^* L_1 \otimes \pi_2^*L_2.
\]
The $\boxtimes$ operation always induces a monomorphism
$\Pic(X_1)\times \Pic(X_2) \hookrightarrow \Pic(X_1\times \Pic(X_2))$. However, typically, not all line bundles
on the product $X_1\times X_2$ are of the form $L_1\boxtimes L_2$. For instance, if $C$ is an elliptic curve, then
the map $\Pic C \times \Pic C \to
\Pic (C\times C)$ is never surjective,
see \cite[\S2.1]{hulek_laface_Pic_numb_of_abelian_vars}.
Similarly, the induced map on the numerical classes
$N^1(X_1)\times N^1(X_2) \hookrightarrow N^1(X_1\times X_2)$ is injective, but not necessarily surjective.

\begin{lemma}
\label{lem_sufficiently_ample_on_product}
Suppose $X_1$ and $X_2$ are projective schemes and $P$ is a property of line bundles on $X_1\times X_2$
which is satisfied by sufficiently ample line bundles.
Then for a sufficiently ample line bundle $L_1$ on $X_1$ and a sufficiently ample line bundle $L_2$ on $X_2$
the property $P$ holds for the line bundle $L_1\boxtimes L_2$ on $X_1\times X_2$.
\end{lemma}

\begin{prf}
Consider the map
$\boxtimes\colon N^1(X_1)\times N^1(X_2) \hookrightarrow N^1(X_1\times X_2)$
of free abelian groups.
Since $X_i$ (and hence also $X_1\times X_2$) are projective, the cone $\Nef(X_1\times X_2)$ is of full dimension.
We apply Lemma~\ref{lem_two_cones_in_lattice_and_sublattice} to
$N=N^1(X_1\times X_2)$,
$N'=N^1(X_1)\times N^1(X_2)$,
$\sigma = \Nef(X_1\times X_2)$,
$\sigma'= \Nef(X_1) \times \Nef(X_2)$
to conclude, that there exists a pair
$(D_1,D_2)\in N^1(X_1)\times N^1(X_2)$ such that $P$ is satisfied for any
$L_1\boxtimes L_2$ with the numerical equivalence class of $L_i$ in
$D_i + \Nef(X_i)$, as claimed.
\end{prf}

\section{Sufficiently ample line bundles}\label{sec_suff_ample_line_bndles}
The main result in the article is about  a property that is fulfilled
for sufficiently ample line bundles.
We commence with a thorough discussion of this notion,
first introduced in~\cite[Def.~3.1]{green-koszul-cohomology-2}, also exploited
in~\cite{sidman_smith_linear_det_eq_for_all_proj_schemes}.

\subsection{Derived push forward and Fujita vanishing}

One of our main tools in the article is the relative Fujita vanishing theorem. In fact, this statement
(Theorem~\ref{thm_relative_fujita_vanishing}) and its classical version (Corollary~\ref{cor_fujita_vanishing})
are the grandparents of many properties that hold for sufficiently ample line bundles.

In order to fully enjoy the involvement of the grandpas we also exploit a result proved in the book of
Mumford \cite[Cor.~2 in Sect.~5, pp.~50--51]{mumford_abelian_varieties}.
Although very useful, its presentation in the book is somewhat hard to digest, as it involves sentences:
``let (\dots) be as above'', and above we find
``let (\dots) be as in the theorem (except (\dots))''.
Thus for our own and reader's convenience (and also for further reference) we decipher the statement and
quote it as Proposition~\ref{prop_mumford_derived_pushforwards}.

For a morphism of schemes $Z\to Y$, a coherent sheaf $\ccF$ on $Z$ and a point $y\in Y$,
we denote  by $\kappa(y)$ the field such that $y=\Spec \kappa(y)$,
and by $\ccF_y = \ccF \otimes_{\ccO_Y} \kappa(y)$ the fibre sheaf.
In the special case where $Z=X\times Y$, the morphism $X\times Y \to Y$ is the projection
and $y=\Spec \kk$ is a closed point, by a slight abuse of notation we think of $\ccF_y$ as of a coherent sheaf on $X$
using the natural isomorphism $X\simeq X\times \set{y}$.

\begin{prop}[{\cite[Cor.~2 in Sect.~5, pp.~50--51]{mumford_abelian_varieties}}]
\label{prop_mumford_derived_pushforwards}
Let $Z$ and $Y$  be a locally Noetherian schemes and let $\ccF$ be a coherent sheaf on $Z$.
Assume $Y$ is reduced and connected.
Then for all integers $q$ the following are equivalent
\begin{enumerate}
\item  \label{item_mumford_constant_function}
the function $y \mapsto \dim_{\kappa(y)} H^q(Z_y, \ccF_y) \ $ is constant,
\item \label{item_mumford_locally_free_and_iso_on_fibres}
$R^q f_*(\ccF) $ is a locally free sheaf
and for all $y \in Y$, the natural map
\[
R^q f_*(\ccF) \otimes_{\ccO_Y} \kappa(y)
\to H^q (Z_y, \ccF_y)
\]
is an isomorphism.
\end{enumerate}
Moreover, if the equivalent conditions \ref{item_mumford_constant_function}
and \ref{item_mumford_locally_free_and_iso_on_fibres} are fulfilled for some $q$,
then also the map
\[
R^{q-1} f_*(\ccF) \otimes_{\ccO_Y} \kappa(y)  \to  H^{q-1}(Z_y, \ccF_y)
\]
is an isomorphism.
\end{prop}

\begin{lemma}
\label{lem_application_of_mumford_to_product}
Suppose $X$  is a projective scheme
and $Y$ be a connected reduced quasiprojective scheme over $\kk$.
Fix a coherent sheaf $\ccF$ on the product $X \times Y$ with $\pi_X\colon X \times Y \to X$
projection to $X$, and $\pi_Y\colon X \times Y \to Y$  projection to $Y$.
Assume $\ccF$ is flat over $Y$ and  such that
$
R^q (\pi_Y)_* \ccF =0
$
for all $q>0$.
Then for any closed point $y=\Spec\kk \in Y$
\[
H^q(X, \ccF_y)=0 \qquad  \text{ for all }  q>0.
\]
\end{lemma}

\begin{prf}
We will use Proposition~\ref{prop_mumford_derived_pushforwards} with $Z=X\times Y$.
Pick a sufficiently large $q$.
Then the condition \ref{item_mumford_constant_function} of Proposition~\ref{prop_mumford_derived_pushforwards}
for that $q$ is satisfied, as the cohomologies vanish for $q$ greater than the dimension of $X$, in particular,
they are constantly equal to $0$.
Hence by the ``moreover'' part of Proposition~\ref{prop_mumford_derived_pushforwards}
the map
\[
R^{q-1}(\pi_Y)_*(\ccF) \otimes \kappa(y) \to H^{q-1}(X,  \ccF_y)
\]
is an isomorphism.
Moreover, provided $q>1$, by our assumptions
$R^{q-1}(\pi_Y)_*\ccF=0$ (in particular, it is locally free),
hence  condition~\ref{item_mumford_locally_free_and_iso_on_fibres}
is satisfied for a smaller value of
$q$ and hence also condition~(i) is satisfied.
By downward induction on $q$, we  get the vanishing $ H^{q}(X,  \ccF_y)=0$
for all $q>0$ as claimed.
\end{prf}

\begin{thm}[Relative Fujita vanishing,
{\cite[Thm~1.5]{keeler-ample-filter}}]\label{thm_relative_fujita_vanishing}
Let $A$ be a commutative Noetherian ring, let $Z$ be a projective scheme
over $A$, with the morphism $\pi\colon Z \to \Spec A $.
For all coherent sheaves $\ccF$ on $Z$,
there exists a $\pi$-ample line bundle $L$ on $Z$ such that
\[
H^q(Z, \ccF \otimes  L \otimes D) =0\]
for $q>0$ and all numerically effective line bundles  $D$.
\end{thm}

The special case, when $A=\kk$, is known as Fujita Vanishing Theorem, see
\cite[Thm~(1)]{fujita_vanishing_theorems} or \cite[Thm~1.4.35]{PAG_I}. It can be also seen as a strengthening
of Serre Vanishing Theorem.
\begin{cor}[Fujita vanishing]\label{cor_fujita_vanishing}
Let $X$ be a projective scheme and fix a coherent sheaf $\ccF$ on $X$.
Then for a sufficiently ample line bundle $L$ on $X$ the higher cohomologies vanish:
\[
H^q(X, \ccF \otimes  L) =0 \text{ for all } q>0.
\]
\end{cor}

Finally, in the notation of Theorem~\ref{thm_relative_fujita_vanishing},
we can interpret the $A$-module $H^q(Z, \ccF \otimes  L \otimes D)$ as the global sections
of the coherent sheaf $R^q\pi_*(\ccF \otimes  L)$ on $\Spec A$.
Thus, for more general projective morphism
$\pi\colon Z\to Y$ and covering $Y$ by open affine subsets, we can conclude
appropriate vanishings for  $R^q\pi_*(\ccF \otimes  L)$ sheaves.
Explicitly, the interesting case is stated in the following lemma.

\begin{lemma}
\label{lem_relative_Fujita_vanishing_for_R_q}
Let $X$  be a projective scheme and $Y$ be a quasiprojective scheme over $\kk$.
Fix a coherent sheaf $\ccF$ on the product $X \times Y$ with $\pi_X$
projection to $X$, and $\pi_Y$  projection to $Y$. Then for a sufficiently ample line bundle $L$ on $X$
\[
R^q (\pi_Y)_*\Bigl( \ccF \otimes \pi_X^*L   \Bigr ) =0 \qquad  \text{ for all }  q>0.
\]
\end{lemma}
\begin{prf}
The statement is local on $Y$, thus we may replace $Y$ with its affine open subsets.
Thus assume $Y=\Spec A$ for some $\kk$-algebra $A$.
Now $L$ is sufficiently ample on $X$ if and only if $\pi_X^* L$ is sufficiently ample on $X\times \Spec A$.
By~\cite[Thm~III.8.5]{hartshorne} we have
\[
R^q (\pi_Y)_*\Bigl( \ccF \otimes \pi_X^*L   \Bigr )(\Spec A) = H^q( X \times \Spec A, \ccF \otimes \pi_X^*L   )
\]
and the latter is $0$ by the relative Fujita vanishing, Theorem~\ref{thm_relative_fujita_vanishing}.
This concludes the proof of the first claim.
\end{prf}

\subsection{Very ampleness}

As a warm up before addressing the issues of cactus varieties we review the fact
that sufficiently ample line bundle is very ample.
Although it is known to some experts, it is not easy to reference explicitly, and perhaps it is worth explaining in more detail,
as it also provides a nice application of relative Fujita vanishing in the version of
Lemma~\ref{lem_relative_Fujita_vanishing_for_R_q}.

\begin{prop}[sufficiently ample line bundles are very ample]\label{prop_sufficiently_ample_is_very_ample}
Let $X$ be a projective variety. Then there exists a very ample line bundle
$L \in \Nef(X)$
such that for any numerically effective line bundle $D \in \Nef(X)$
the line bundle  $L\otimes D$ is very ample.
\end{prop}

Before discussing the proof we comment a little on the content. In general, if $L_1$ is very ample and $L_2$ is nef
or even ample, then $L_1\otimes L_2$ is not necessarily very ample, see Example~\ref{ex_hypereliptic_curve}.
Moreover, the same example illustrates that \emph{very ampleness} is not invariant under numerical equivalence,
that is if $L_1$ is very ample and $L_2\equiv L_1$,
then $L_2$ is not necessarily very ample.
However, the proposition shows that these two pathologies can only happen ``close to the boundary'' of the nef cone.
Instead, if we are sufficiently deep in the interior of $\Nef(X)$, then all line bundles are very ample,
as illustrated by the following figure.
\[
\begin{tikzpicture}[baseline=.4cm,scale=0.5]
\fill[color=czerwony!20] (0,0) -- (4.5,4.5) -- (9.5, 4.5)--(9.5,0)
-- cycle;
\draw[step=1, ciemnyblekit, thin] (-2.5,-0.5) grid (9.5, 4.5);
\draw[line width=0.5mm, czerwony, densely dotted] (4.5,4.5)-- (0,0)  -- (9.5,0)  (2.5,1) node[anchor=north west]{ample cone};
\fill[thick, granat] (0,0) node[anchor=south east]{$\ccO_X$} circle (6pt);
\fill[color=zielony!20] (3,2) -- (5.5,4.5) -- (9.5, 4.5)--(9.5,2)
-- cycle;
\draw[line width=0.5mm, zielony] (5.5,4.5)-- (3,2)  -- (9.5,2)  (3.5,3.2) node[anchor=north west]{very ample!};
\fill[thick, zielony] (2,1) circle (6pt);
\fill[thick, zielony] (5,1) circle (6pt);
\end{tikzpicture}
\] 

\begin{example}\label{ex_hypereliptic_curve}
Suppose $\cchar \kk \neq 2, 3$,
and consider a smooth plane curve
\[
C = \set{x^4+ y^4 + yz^3 =0}\subset \PP^2
\]
of genus $g=3$,
and let $D = [0,0,1] \in C$ be a $\kk$-point of $C$ seen as a prime divisor on $C$.
Then the hyperplane section $(y=0)$ is the divisor $4D$ on $C$,
thus the corresponding line bundle $\ccO_C(4D)$
is very ample and $\ccO_C(D)$ is ample.
Moreover, $\ccO_C(4D)=\ccO_{\PP^2}(1)|_{C} = (\ccO_{\PP^2}(-3) \otimes \ccO_{\PP^2}(C))|_{C} $ and  by the adjunction formula $4D$ is the canonical divisor  of $C$.
On the other hand, $\ccO_C(5D)$ is not very ample.
Indeed, by Riemann-Roch \cite[Thm~IV.1.3]{hartshorne} we have
\begin{align*}
h^0(\ccO_C(4D)) & = \deg (4D)+1 - g + h^0(\ccO_C) = 4+1-3+1= 3,\\
h^0(\ccO_C(5D)) & = \deg (5D)+1 - g + h^0(\ccO_C(-D)) = 5+1-3+0= 3.
\end{align*}
If $s\in H^0(\ccO_C(D))$ is the non-zero section defining $D$,
then the multiplication by $s$ determines an injective map
$H^0(\ccO_C(4D)) \stackrel{\cdot s}{\hookrightarrow} H^0(\ccO_C(5D))$.
By the dimension count above, this map is an isomorphism, thus every section of
$\ccO_C(5D)$ vanishes at $D$ and thus $\ccO_C(5D)$ cannot be very ample.
Moreover, if $L$ is another line bundle of degree $4$, not isomorphic to $\ccO(4P)$,
then $L$ is also not very ample, since $h^0(L)=2$.

The following table shows which divisors of the form $kD$ on $C$ are very ample (\emph{v.a.}), and which are not.
Moreover, any divisor of degree at least $7$ is very ample, illustrating the statement of
Proposition~\ref{prop_sufficiently_ample_is_very_ample} for the case of $X=C$.

\medskip

\noindent
\begin{tabular}{c|c|c|c|c|c|c}
  $\ccO_C(D)$ & $\ccO_C(2D)$ & $\ccO_C(3D)$& $\ccO_C(4D)$& $\ccO_C(5D)$& $\ccO_C(6D)$ & $\ccO_C(E),
                                \deg E \geqslant 7$  \\
\hline
  \textcolor{czerwony}{not v.a.}& \textcolor{czerwony}{not v.a.}&\textcolor{czerwony}{not v.a.}& \textbf{v.a.}
  &\textcolor{czerwony}{not v.a.}& \textcolor{czerwony}{not v.a.} & \textbf{v.a.}
\end{tabular}

\medskip
\end{example}

We now relate very ampleness to vanishing of certain $H^1$-cohomologies.
For any subscheme $R\subset X$ denote by $\ccI_R\subset \ccO_X$ the ideal sheaf of $R$.

\begin{lemma}\label{lem_criterion_for_very_ampleness_all_schemes}
Suppose $X$ is a projective scheme and $L$ is a line bundle on $X$.
If for any finite subscheme $R\subset X$ of degree $2$ we have  $H^1(X, \ccI_R\otimes L) =0$, then $L$ is very ample.
\end{lemma}

\begin{prf}
A line bundle $L$ is very ample when it separates points and tangent directions.
More precisely, when for any $0$-dimensional scheme $R$ of degree $2$ (which is either isomorphic
to two disjoint reduced points $\Spec( \kk\times \kk)$ or to the double point on a line
$\Spec \kk[\varepsilon]/(\varepsilon^2)$), the restriction map
$
H^0(X,L) \to H^0(R, L|_R)
$
is surjective.
Considering the twisted ideal-ring short exact sequence
\[
0 \to \ccI_R\otimes L \to L \to L|_R \to 0
\]
we obtain the long exact sequence of cohomology
\[
0 \to H^0(X, \ccI_R\otimes L)\to H^0(X, L)) \to H^0(X,L|_{R})\to H^1(X, \ccI_R\otimes L)\to \dotsc
\]
Since the $H^1$ term vanishes by the assumption of the Lemma,
it follows that the second $H^0$ map is surjective, as desired.
\end{prf}

By Fujita vanishing (Corollary~\ref{cor_fujita_vanishing}) for each scheme $R$ we have
the vanishing of $H^1(I_R\otimes L)$ for sufficiently ample $L$.
Thus for each $R$ we obtain a translated cone $\sigma_R = D_R + \Nef(X)$ (for some divisor class $D_R$).
But at this point it is not clear that the intersection $\bigcap_{R} \sigma_R$ is non-empty, as there are infinitely
many $R$ as soon as $\dim X > 0$.
Thus we need some way of controlling the vanishing all at the same time, and this is where
the relative version of Fujita vanishing comes in handy.

Let $\ccH:=\reduced{(\usualHilb_2 X)}$ be the reduced subscheme of the Hilbert scheme of degree $2$
subschemes of $X$, so that its closed points
$[R]\in \ccH$ are in 1-1 correspondence to degree $2$ subschemes $R\subset X$ finite over $\kk$.
Consider the product $X \times \ccH$ together with
two natural projections $\pi_X\colon X \times \ccH \to X$ and $\pi_{\ccH}\colon X \times \ccH \to \ccH$.
The \emph{universal ideal} $\ccJ\subset \ccO_{X \times \ccH}$ is a flat over $\ccH$ ideal sheaf,
whose fibre $\ccJ_{[R]}\subset \ccO_X$ is the ideal sheaf $\ccI_{R} \subset \ccO_X$.

\begin{lemma}
\label{lem_criterion_for_very_ampleness_relative_push_forward}
Suppose $X$ is a projective scheme, $L$ is a line bundle on $X$, and
$\ccH=\reduced{(\usualHilb_2 X)}$ and $\ccJ$ is as above.
If $R^q (\pi_{\ccH})_*(\ccJ \otimes \pi_X^*L) =0$ for all $q> 0$,
then $L$ is very ample.
\end{lemma}

\begin{prf}
  A line bundle is very ample if and only if it is very ample on each of its connected components.
  Thus without loss of generality we may assume $X$ is connected, and therefore also $\ccH$ is connected.
By definition $\ccH$ is also reduced and projective.
The ideal sheaf $\ccJ$ is flat by definition of the Hilbert scheme and of the universal ideal.

We apply Lemma~\ref{lem_application_of_mumford_to_product}
with $Y=\ccH$ and $\ccF =\ccJ\otimes \pi_X^*L$.
Thus $H^q(X, \ccJ_{[R]}\otimes L) = 0$ for all $q>0$ and all $[R]\in \ccH$.
The definition of the universal ideal $\ccJ_{[R]} = \ccI_R$ and the vanishing of the $H^1$-cohomologies imply
that $L$ is very ample by Lemma~\ref{lem_criterion_for_very_ampleness_all_schemes}.
\end{prf}

\begin{prf}[ of Proposition~\ref{prop_sufficiently_ample_is_very_ample}]
  Let $\ccH = \reduced{(\usualHilb_2 X)}$, $\ccJ$ be the universal ideal sheaf on $X\times \ccH$, and
  the projections $\pi_X$, $\pi_{\ccH}$ be as above.
We apply  Lemma~\ref{lem_relative_Fujita_vanishing_for_R_q}
with $Y=\ccH$ and $\ccF=\ccJ$.
Thus for a sufficiently ample line bundle $L$ on $X$ we have
\[
R^q (\pi_{\ccH})_*(\ccJ \otimes \pi_X^*L) =0 \text{ for all } q> 0,
\]
and by Lemma~\ref{lem_criterion_for_very_ampleness_relative_push_forward}
we get that $L$ is very ample.
\end{prf}

For future reference we state and prove another lemma on surjectivity of multiplication of section.

\begin{lemma}
\label{lem_multiplication_of_sections_is_surjective_for_sufficiently_ample}
Suppose $X$ is a projective scheme.
The multiplication of sections
\[
H^0(X,L_1)\otimes H^0(X,L_2) \to H^0(X,L_1\otimes L_2)
\]
is surjective
when  $L_1$ and $L_2$ are  sufficiently ample line bundles.
\end{lemma}

\begin{prf}
Consider the product $X \times X$     and the diagonal embedding
$\Delta \colon X\to X \times X$.  Then for any line bundles $L_1$ and $L_2$ on $X$,     we have
\begin{equation}
\label{equ_isomorphisms_on_diagonal}
(L_1 \boxtimes L_2)\otimes_{\ccO_{X\times X}}  \ccO_{X} = \Delta^* (L_1 \boxtimes L_2) =  L_1\otimes L_2.
\end{equation}
Let $\ccI_{\Delta}\subset \ccO_{X\times X}$ be the ideal sheaf of $\Delta(X)$.
We have a short exact sequence:
\[
0 \to \ccI_{\Delta}
\to  \ccO_{X\times X} \to \ccO_X \to 0.
\]
Twisting it by $L_1 \boxtimes L_2$ and exploiting \eqref{equ_isomorphisms_on_diagonal} we obtain:
\[
0 \to \ccI_{\Delta}\otimes (L_1 \boxtimes L_2)
\to  L_1 \boxtimes L_2 \to L_1 \otimes L_2 \to 0
\]
Taking into account that
$H^0(X \times X, L_1 \boxtimes L_2) = H^0(X,L_1)\otimes H^0(X,L_2)$,
the corresponding long exact sequence of cohomologies includes:
\begin{equation}
\label{equ_cohomology_exact_sequence_for_diagonal}
\dotsb \to
H^0(X,L_1)\otimes H^0(X,L_2) \to H^0(L_1 \otimes L_2) \to H^1\left(\ccI_{\Delta}\otimes (L_1 \boxtimes L_2)\right) \to \dotsb
\end{equation}
where the first map is the multiplication of sections by construction.
Therefore, whenever $H^1\left(\ccI_{\Delta}\otimes (L_1 \boxtimes L_2)\right) =0$,  then
the multiplication of sections is surjective.

By Fujita vanishing (Corollary~\ref{cor_fujita_vanishing}),
for a sufficiently ample line bundle $L$ on $X\times X$ we have
$H^1\left(\ccI_{\Delta}\otimes L\right) =0$.
Moreover, by Lemma~\ref{lem_sufficiently_ample_on_product},
we can take $L=L_1 \boxtimes L_2$ for sufficiently ample line bundles $L_1$ and $L_2$ on $X$.
Thus $H^1\left(\ccI_{\Delta}\otimes (L_1 \boxtimes L_2)\right) =0$,
and the claim of the lemma follows from the exact sequence
\eqref{equ_cohomology_exact_sequence_for_diagonal}.
\end{prf}

\section{Multigraded rings of sections}
\label{sec_multigraded_rings}
Let $X\subset \PP^N$ be an embedded projective scheme.
We will discuss divisorial algebras of sequences of divisors on $X$,
similarly to \cite[\S1.3.2]{arzhantsev_derenthal_hausen_laface_Cox_rings} or
\cite[\S1]{kuronya_urbinati_geometry_of_multigraded_rings_and_embeddings}, which however have somewhat
different settings (among other differences, these two references require $X$ to be a  normal variety or prevariety).
The only cases we will be interested in here are those of a single divisor and those of a sequence of two divisors.
In the first case, the traditional name is the \emph{section ring}, and we will stick to this one instead of
the more general ``divisorial algebra''. In the second case we use \emph{double section ring} for internal consistency.
Moreover, the only cases we work with here are the cases of very ample Cartier divisors.
By a slight abuse of notation, we will work with the corresponding line bundles instead of divisor,
being aware that strictly speaking the $\kk$-algebra structures
(but not their roles in this paper or isomorphism types) depend on
the uniform choices of the divisors in the appropriate linear systems.
We now discuss in detail the two cases (single or double) separately.

\subsection{A variant of projective normality}

For any line bundle $L$ on $X$ denote by $S_L$
the \emph{section ring} $S_L:= \kk\oplus \bigoplus_{d> 0} H^0(X, L^{\otimes d})$, which is an $\NN$-graded $\kk$-algebra%
\footnote{It might be somewhat controversial and/or unorthodox to define the zeroth grading as $\kk$ instead of
$H^0(X, \ccO_X)$. If $X$ is connected and reduced, then $H^0(X, \ccO_X) =\kk$ and the controversy does not appear.
Otherwise, $H^0(X, \ccO_X)$ might be a different finite $\kk$-algebra.
However, for consistency issues
it is more convenient for us to use just $\kk$.}.
Let us fix an embedding of $X$ via the 
complete linear system of a line bundle into projective space $\PP^N$. We denote
$S[\PP^N] = \kk[\fromto{\alpha_0}{\alpha_N}]$ be the standard graded polynomial ring and at the same time
the \emph{homogeneous coordinate ring} of $\PP^N$.
Let  $I_X \subset S[\PP^N]$  be the homogeneous ideal of $X$ in this embedding.
Thus, when $L$ is the restriction of $\ccO_{\PP^N}(1)$ to $X$, we have two graded rings:
\begin{itemize}
\item  the classical homogeneous coordinate rings of $X$, denoted $S[X]: = S[\PP^N]/I_X$, and
\item the section ring $S_L$ as above.
\end{itemize}
We always have an embedding $S[X] \subset S_L$.

Moreover, $X$ is called \emph{projectively normal} if the affine cone $\Spec S[X]$ is a normal variety.
Equivalently, by \cite[Ex.~II.5.14(d)]{hartshorne},
$X$ is projectively normal if and only if
\begin{itemize}
\item $X$ is a normal variety and
\item the inclusion above is an equality $S[X] = S_L$.
\end{itemize}

Here we are more interested in the second condition, $S[X] = S_L$ for more general schemes than normal varieties.
Occasionally, somewhat confusingly, in informal conversations schemes satisfying $S[X] = S_L$ are also called
projectively normal, but we refrain from using this name as it has nothing to do with the structure of $X$ itself,
in particular, nothing to do with normality of $X$.

Since we are primarily interested in complete embeddings of $X$ given by $H^0(X,L)$ for a very ample $L$,
we will define this as a property of $L$.

\begin{defin}
We say a line bundle $L$ on a projective scheme $X$ is \emph{standard graded}
if $L$ is very ample and the section ring $S_L$ as a graded $\kk$-algebra is generated in the first degree,
that is $H^0(X,L)$ generates $S_L$.
\end{defin}
The name ``standard graded'' refers to the traditional name ``standard graded $\kk$-algebra'',
see for instance \cite[p.~115]{stanley_combinatorics_and_commutative_algebra},
\cite[\S2.2]{jelisiejew_PhD}.
Thus the full name for this notion should be:
``$L$ is very ample and has standard graded ring of sections'', which we abbreviate to
simply ``$L$ is standard graded''.

As discussed above, whenever $L$ is standard graded, we have a presentation of $S_L$ as
a quotient $S_L= \kk[\fromto{\alpha_0}{\alpha_N}]  \slash I_X $ 
for the homogeneous ideal $I_X \subset \kk[\fromto{\alpha_0}{\alpha_N}]$, which is
the homogeneous saturated ideal of $X$ under the embedding $X\hookrightarrow \PP^N = \PP(H^0(X,L)^*)$.

Here the grading of the polynomial ring is standard, that is $\deg \alpha_j=1$ for each $j\in \setfromto{0}{N}$,
the number of variables $N+1$ is equal to the dimension of the linear system $H^0(X,L)$.
The reason we stress the above fact is that thanks to this presentation we are allowed to exploit Hilbert function
tricks such as Macaulay bound on the growth for Hilbert function and Gotzmann persistence of the growth
for both $S_L$ and its quotients by homogeneous ideals, as illustrated in the  standard
Lemma~\ref{lem_Gotzmann_in_single_graded}.

\begin{defin}\label{def_d_embedding}
Suppose $X$ is a projective scheme and  $L$ is standard graded line bundle. We say $L$
\emph{is a  $d$-embedding}
if  the ideal $I_X$ as above is generated by elements of degree  at most $d$.
\end{defin}

For a homogeneous ideal $I\subset S_L$ the subscheme in $X$ defined by $I$ is denoted by $Z(I)$.
\begin{lemma}
\label{lem_Gotzmann_in_single_graded}
Suppose $I\subset S_L$ is a homogeneous ideal such that for some $r,d$ with $d\geqslant r$
the following assumptions hold:
\begin{itemize}
\item $L$ is a $d$-embedding
\item $\dim_{\kk} (S_L/I)_d =\dim_{\kk} (S_L/I)_{d+1} = r$, and
\item $I$ is generated by elements of degree at most $d$.
\end{itemize}
Then:
\begin{itemize}
\item $Z(I)$ is a finite subscheme of $X$ of degree  $r$, and
\item the linear span of $Z(I) \subset X\subset \PP\left(H^0(X,L^{\otimes d'})^*\right)$ is equal
to $\PP\left(I_{d'}^{\perp}\right)$ for all $d'\geqslant d$.
\item $I$ agrees with its saturation from degrees $d$ onwards: $I_{d'}  = I(Z(I))_{d'}$ for $d'\geqslant d$.
\end{itemize}
\end{lemma}

\begin{prf}
Let  $\Phi \colon \Bbbk[x_0,\ldots,x_n] \to S_L $ be the quotient morphism whose kernel is the ideal $I_X$,
so that $I_X=\Phi^{-1}((0))$. In the polynomial ring we consider  the grading by the degree of a polynomial
and the induced grading in  the quotient ring $S_L$.
Let $\widehat{I} = \Phi^{-1}(I)$ be the preimage of the ideal $I \subset S_L$.
Since $L$ is a $d$-embedding, the ideal $\widehat I$ is generated by elements of degree at most $d$.
Also, by our second assumption on $I$ we have
\[
\dim_{\kk} (\kk[x_0,\ldots,x_N]/\widehat {I})_d =\dim_{\kk} (\kk[x_0,\ldots,x_N]/\widehat I)_{d+1} = r
\]
which means the Hilbert function of the scheme $Z(\widehat{I})$ satisfies
\[
H_{Z(\widehat{I})}(d)=H_{Z(\widehat{I})}(d+1)=r.
\]
Since  $d \geqslant r$, the assumptions of the Gotzmann persistence theorem, 
see~\cite[Thm~4.3.3]{bruns_herzog_Cohen_Macaulay_rings}, are satisfied
--- the maximal jump here is zero  and it is attained, so it persists and
\[
\dim_{\kk} (\kk[x_0,\ldots,x_N]/\widehat {I})_{d'} =\dim_{\kk} (S_L/I)_{d'} = r \text{ for all } d' \geqslant d.
\]
The Hilbert function of $Z(I)$ is the same, so we proved the first part of the claim.
Finally, the  linear span in the embedding given by $L^{\otimes d'}$ is the zero set
of the elements of degree $d'$ that vanish on $Z(I)$ and its dimension is equal to $H_{Z(I)}(d')$.
\end{prf}

An expert reader is probably aware that sufficiently ample line bundles are standard graded.
This statement, Corollary~\ref{cor_standard_graded_for_suff_ample}\ref{item_single_standard_graded},
serves us merely as an example, and we will prove it later as a consequence of a stronger result,
which is relevant in our context.

\subsection{Double section ring}

In order to prove our results in Section~\ref{sec_cacti} we also need to consider \emph{double section ring}.
That is,   for two line bundles $L_1$ and $L_2$ we consider the bigraded algebra:
\[
S_{L_1,L_2}:= \kk \oplus
\bigoplus_{(d,e)\in \NN^2, \ (d,e)\neq (0,0)} H^0(X, L_1^{d}\otimes L_2^e).
\]
Here, similarly as in the definition of $S_L$, in the $(0,0)$-grading we place $\kk$
instead of seemingly more natural $H^0(X,\ccO_X)$.
Moreover, as mentioned above, in order to define the ring structure on $S_{L_1, L_2}$
we need to make a choice of divisors $H_i$ in the linear system of $L_i$ and consistently, a choice of  $d H_1 + e H_2$
in the linear system of $L_1^{d}\otimes L_2^e$.
Since different consistent choices lead to isomorphic $\kk$-algebras,  we skip $H_1$ and $H_2$ from our notation.

\begin{defin}
\label{def_doubly_standard_graded}
For a projective scheme $X$ and a pair of line bundles $L_1, L_2$ on $X$ we say that $(L_1, L_2)$
is \emph{doubly standard graded} if
both $L_1$ and $L_2$ are very ample and $S_{L_1,L_2}$ is generated by $(1,0)$ and $(0,1)$ gradings.
\end{defin}

Note that $S_{L_1}\subset S_{L_1, L_2}$, as the sum of $(d, 0)$-degree parts  for all $d\in \NN$.
Thus if $(L_1, L_2)$ is doubly standard graded, then $L_1$ (and similarly $L_2$) is standard graded.
In fact, any $L_1^{\otimes d}\otimes L_2^{\otimes e}$ is standard graded for any $d,e\in \NN^2\setminus\set{(0,0)}$.

\begin{lemma}
\label{lem_doubly_standard_graded_algebra}
Suppose $S = \bigoplus_{d,e \geqslant 0}$ is a bigraded algebra with $S_{(0,0)} = \kk$ and
generated by $S_{(1,0)}$ and $S_{(0,1)}$.
Then for all $d\geqslant 0$ and $e\geqslant 0$ the multiplication
maps $S_{(d,e)} \otimes S_{(1,0)} \to S_{(d+1,e)}$ and $S_{(d,e)} \otimes S_{(0,1)} \to S_{(d,e+1)}$ are surjective.
\end{lemma}

\begin{prf}
It is enough to prove surjectivity of one map, the other follows by swapping the components of the grading.
We have the surjective map $S_{(1,0)}^{\otimes (d+1)}\otimes S_{(0,1)}^{\otimes e} \to  S_{(d+1,e)}$,
which factors through
\[
S_{(1,0)}^{\otimes (d+1)}\otimes S_{(0,1)}^{\otimes e} \to S_{(d,e)}\otimes S_{(1,0)} \to  S_{(d+1,e)}
\]
forcing also the latter map to be surjective, as claimed.
\end{prf}

By analogy to the (single) standard graded case,
if $(L_1, L_2)$ is doubly standard graded,
then we have a presentation of the double section ring $S_{L_1, L_2}$ as a quotient of a bigraded polynomial ring:
\begin{equation}
S_{L_1, L_2} = \kk[\fromto{\alpha_0}{\alpha_{N_1}}, \fromto{\beta_0}{\beta_{N_2}}]/J_{X},
\end{equation}
where $\dim H^0(X, L_i) = N_i$
for $i=1,2$,
$\deg \alpha_j = (1,0)\in \ZZ^2$ for $j = \fromto{0}{N_1}$,
$\deg \beta_k = (0,1)\in \ZZ^2$ for $k = \fromto{0}{N_2}$,
$J_X$ is a bihomogeneous ideal. Here $J_X$ can be seen as the ideal of $X$
embedded in $\PP^{N_1}\times \PP^{N_2}$ using both linear systems of $L_i$.
If $R\subset X$ is a subscheme,
then define the bihomogeneous ideal $J_R \subset S_{L_1, L_2}$ to be generated by
all sections $s\in H^0(X, L_1^{d}\otimes L_2^e)$ --- for any $(d,e)$ --- that vanish identically on $R$.
\begin{lemma}
\label{lem_ideal_of_scheme_is_saturated}
For any closed subscheme $R \subset X$ the ideal $J_R$ is saturated
with respect to both  $\left( S_{(1,0)} \right)$ and $\left( S_{(0,1)} \right)$.
\end{lemma}
\begin{prf}
$J_R=\set{ f \in S_{L_1,L_2} \mid  f_{\mid R} =0 }$ and its saturation with respect to $\left( (S_{L_1,L_2})_{(1,0)} \right)$
is the set $\set{f \in S_{L_1,L_2} \mid \forall_i \exists_k \quad f \cdot \alpha_i^k\in J_R}$.
Since $S_{L_1,L_2}$ is generated by gradations $(1,0)$ and $(0,1)$,
the scheme $X$ is covered by the affine pieces $X = \bigcup \ccU_{i,j}$, where
$\ccU_{i,j} $ are the sets  $\set{x \in X \mid \alpha_i(x) \neq 0 \text{ and } \beta_j(x) \neq 0 }$.
Now take an  $f \in \left((S_{L_1,L_2})_{d,e}\right)$ with $f \in \left(J_R : \left((S_{L_1,L_2})_{(1,0)}\right)^{\infty} \right)$.
This means for some number $k$ and all $i \in \set{1,\ldots,N_1}$ we have $f \cdot \alpha_i^k \in J_R$
where $\alpha_i$ is invertible on $\ccU_{i,j}$.
Since by definition $J_R= \bigoplus_{d,e}H^0(\ccI_R \otimes \ccO_X(d,e))$ and we have proved that $f$
vanishes on each $\ccU_{i.j} \cap R$, we obtain $f \in J_R$.
\end{prf}
\begin{lemma}\label{lem_Gotzmann_and_Hilbert_function}
Suppose $(L_1,L_2)$ is doubly standard graded on $X$
and  $R\subset X$ is a finite subscheme of degree $r$.
Then for any $(d,e)$:
\begin{enumerate}
\item
\label{item_lenght}
$\dim (S_{L_1,L_2}/J_R)_{d,e} \leqslant r$,
\item
\label{item_hilbert_polynomial}
if $d+e\ge r-1$ then $\dim (S_{L_1,L_2}/J_R)_{d,e} =r$,
\item
\label{item_monotone_nonzerodivisor_exists}
$\dim (S_{L_1,L_2}/J_R)_{d,e} \leqslant \dim (S_{L_1,L_2}/J_R)_{d+1,e}$,
$\dim (S_{L_1,L_2}/J_R)_{d,e} \leqslant \dim (S_{L_1,L_2}/J_R)_{d,e+1}$,
\item
\label{item_version_of_Gotzmann}
if $\dim (S_{L_1,L_2}/J_R)_{d,e} = \dim (S_{L_1,L_2}/J_R)_{d+1,e}$
then
\[
\dim (S_{L_1,L_2}/J_R)_{d',e} = \dim (S_{L_1,L_2}/J_R)_{d,e} \text{ for any } d'\geqslant d,
\]
\item
\label{item_equations_of_linear_span}
for $X \subset \PP(H^0(X, L_1^{d}\otimes L_2^e)^*)$, the linear span
$\linspan{R}$ in this projective space
is defined by $(J_R)_{d,e}$.
\end{enumerate}
\end{lemma}
\begin{prf}
For the  item~\ref{item_lenght} recall  $\ccI_R$ is the ideal sheaf of $R$. We have the exact sequence
\begin{equation}\label{equ_exact_seq_double_graded}
0 \to H^0(  \ccI_R \otimes L_1^{d} \otimes L_2^{e}) \to H^0(L_1^d \otimes L_2^e) \to
H^0( \ccO_R \otimes L_1^d \otimes L_2^e) \to H^1(  \ccI_R \otimes L_1^{d} \otimes L_2^{e})      
\end{equation}
the beginning of which is the same as
$    0 \to \left( J_R \right)_{(d,e)} \to S_{(d,e)} \to  H^0(\ccO_R)$.
For the identification of fourth entries we pick an isomorphism
$  \ccO_R \otimes L_1^d \otimes L_2^e \simeq \ccO_R$.
Since $R$ is a finite scheme, $\dim H^0( \ccO_R) =r$.

Now we prove~\ref{item_monotone_nonzerodivisor_exists}. The statement is implied by the existence
of non-zero divisors in the quotient ring $S_{L_1,L_2} / J_R$ in the gradation $(1,0)$ or $(0,1)$ respectively.
Consider primary decomposition of the ideal $J_R= \mathfrak p_1 \cap \ldots \cap \mathfrak p_s$.
By~\cite[Lemma 00LD]{stacks_project}, the set of zero divisors of $S_{L_1,L_2}/J_R$ is the finite union
$\bigcup \left(\sqrt {\mathfrak p_i}/J_R\right)$. We only need to know that
$\bigcup (\sqrt{\mathfrak p_i})_{(1,0)} \subsetneq S_{(1,0)}$.
By Lemma~\ref{lem_ideal_of_scheme_is_saturated} for each $i \in \set{1,\ldots,s}$ it is true that
$(\sqrt{\mathfrak p_i})_{(1,0)} \subsetneq S_{(1,0)}$, so as our field is infinite, the finite union
of vector subspaces of smaller dimension cannot fill a vector space of  bigger dimension.

To prove~\ref{item_version_of_Gotzmann}, observe that  if $f \in (S_{L_1,L_2}/J_R)_{(1,0)}$
is not a zero divisor and the dimensions are equal, the  monomorphism
$(S_{L_1,L_2}/J_R)_{(d,e)} \stackrel{\cdot f}{\to} (S_{L_1,L_2}/J_R)_{(d+1,e)}$  becomes an isomorphism.
Moreover, multiplying by this $f$ in the next gradation
$(S_{L_1,L_2}/J_R)_{(d+1,e)} \stackrel{\cdot f}{\to} (S_{L_1,L_2}/J_R)_{(d+2,e)}$ is an epimorphism,
by a direct computation. Indeed, take any $g \in (S_{L_1,L_2}/J_R)_{(d+2,e)} $ and write $g = \sum \alpha_i \cdot \theta_i$
where each $ \theta_i \in (S_{L_1,L_2}/J_R)_{(d+1,e)}$ can be written as $f\cdot \phi_i$ with $\phi_i\in (S_{L_1,L_2}/J_R)_{(d,e)}$. So we have
$g =f \cdot  \sum \alpha_i \cdot \phi_i$ as claimed,
and multiplying by $f$ is an epimorphism,
but it is also a monomorphism, as $f$ is not a zero-divisor.
Thus the multiplication by $f$ is an isomorphism $(S_{L_1,L_2}/J_R)_{(d+1,e)} \to (S_{L_1,L_2}/J_R)_{(d+2,e)}$.

By Fujita vanishing we know that the last entry of the exact sequence~(\ref{equ_exact_seq_double_graded})
is zero for $d$ or $e$ big enough.
This together with items~\ref{item_version_of_Gotzmann} and~\ref{item_monotone_nonzerodivisor_exists}
implies~\ref{item_hilbert_polynomial}.
The item~\ref{item_equations_of_linear_span} says that the linear
--- from the point of view
of the projective space on the sections of  $L_1^d\otimes L_2^e$  ---
functions that vanish on $R$
are the sections of this line bundle in the ideal of $R$.
\end{prf}

\begin{cor}
\label{cor_standard_graded_for_suff_ample}
Suppose $X$ is a projective scheme.
Then:
\begin{enumerate}
\item \label{item_single_standard_graded}
a sufficiently ample line bundle $L$
is standard graded,
\item \label{item_doubly_standard_graded}
for sufficiently ample line bundles
$L_1, L_2$ the pair $(L_1, L_2)$ is doubly standard graded.
\end{enumerate}
In particular, if $X$ is in addition
a normal variety, then
(for a sufficiently ample $L$)
the embedding of $X$ by a complete linear
system $H^0(X,L)$ is projectively normal.
\end{cor}

In a groaner joke, a mathematician is given the task of removing two nails.
One of the nails is completely hammered into a wooden wall, and the other is stuck out half way.
The mathematician works for hours on an elegant way to pull out the nail that's nailed in.
After finally removing the difficult one, he hammers in completely the one sticking out in order to reduce
to the previous problem.
This is exactly what we are going to do with item \ref{item_single_standard_graded} of the corollary.

\begin{prf}
Item \ref{item_single_standard_graded} follows from \ref{item_doubly_standard_graded}.
The ``in particular'' part follows from
\ref{item_single_standard_graded} and \cite[Ex.~II.5.14(d)]{hartshorne}.

Pick $D\in N^1(X)$ such that for all line bundles $L_1$ and  $L_2$ on $X$ that have their
numerical classes in $D+\Nef(X)$ both of the following properties hold:
\begin{itemize}
\item $L_1$ and $L_2$ are very ample, and
\item the multiplication of sections $H^0(X,L_1)\otimes H^0(X,L_2) \to H^0(X,L_1\otimes L_2)$ is surjective.
\end{itemize}
Such $D$ exists by Lemma~\ref{lem_finitely_many_properties_of_sufficiently_ample_lbs},
Proposition~\ref{prop_sufficiently_ample_is_very_ample}, and
Lemma~\ref{lem_multiplication_of_sections_is_surjective_for_sufficiently_ample}.

To prove \ref{item_doubly_standard_graded},
pick any $L_1$ and $L_2$
such that their numerical classes
are in $D+\Nef(X)$.
Then for any integers $d, e\ge 0$ such that  $(d,e)\ne (0,0)$,
also the class of
$L_1^{\otimes d}\otimes L_2^{\otimes e}$ is in
$D+\Nef(X)$, as it is equal to either:
\begin{itemize}
\item $L_1 \otimes (L_1^{\otimes (d-1)}\otimes L_2^{\otimes e})$ with $d>0$, so $L_1$ is in the translated cone
$D+\Nef(X)$ and $(L_1^{\otimes (d-1)}\otimes L_2^{\otimes e})$ is nef, or
\item $L_2 \otimes (L_1^{\otimes d}\otimes L_2^{\otimes (e-1)})$ with $e>0$, so $L_2$ is in the translated cone
$D+\Nef(X)$ and $(L_1^{\otimes d}\otimes L_2^{\otimes (e-1)})$ is nef.
\end{itemize}

Thus $L_1$ and $L_2$ are very ample
and the multiplication maps
\begin{align*}
H^0(X,L_1^{\otimes (d-1)} \otimes L_2^{\otimes e})\otimes H^0(X,L_1)
& \to H^0(X,L_1^{\otimes d}\otimes L_2^{\otimes e})  \text{ for } d\geqslant 1 , \text{ and}\\
H^0(X,L_1^{\otimes d} \otimes L_2^{\otimes (e-1)})\otimes H^0(X,L_2)
& \to H^0(X,L_1^{\otimes d}\otimes L_2^{\otimes e}) \text{ for } e\geqslant 1 ,
\end{align*}
are surjective, and therefore, the $(d,e)$-th degree of the algebra $S_{L_1,L_2}$
is generated by lower degrees (whenever $d+e\geqslant 2$).
This implies that $S_{L_1,L_2}$ is generated by $H^0(X,L_1)$ and $H^0(X, L_2)$, as claimed
in \ref{item_doubly_standard_graded}.
\end{prf}
\subsection{Annihilators and apolarity in the multigraded setting}\label{sec_apolarity}
Following~\cite{galazka_multigraded_apolarity},
we recall the apolarity tools in the multigraded setting.
There the variety for which the apolarity is developed is a projective toric variety over $\CC$,
instead we need it only for $\PP^{N_1} \times \PP^{N_2}$, but over any algebraically closed $\kk$.
We thus consider the $\ZZ^2$-graded ring $\Bbbk[\alpha_0,\ldots,\alpha_{N_1},\beta_0,\ldots,\beta_{N_2}]$ as before
together with its dual ring, which is the divided power ring in variables $x_0,\ldots, x_{N_1}, y_0, \ldots, y_{N_2}$.
The  apolarity pairing is defined by the following properties:
\begin{multline*}
  \hook \colon \Bbbk[\alpha_0,\ldots,\alpha_{N_1},\beta_0,\ldots,\beta_{N_2}] \times
  \Bbbk_{dp}[x_0,\ldots, x_{N_1}, y_0, \ldots, y_{N_2}] \\
\to  \Bbbk_{dp}[x_0,\ldots, x_{N_1}, y_0, \ldots, y_{N_2}],
\end{multline*}
\begin{align*}
\alpha_i \hook \left(x_{0}^{(a_0)} \dotsm x_{N_1}^{(a_{N_1})} \cdot y_{0}^{(b_0)} \dotsm y_{N_2}^{(b_{N_2})}\right) & =
\begin{cases}
x_{0}^{(a_0)} \dotsm x_i^{(a_i-1)} \dotsm x_{N_1}^{(a_{N_1})} \cdot y_{0}^{(b_0)} \dotsm y_{N_2}^{(b_{N_2})} & \text{if } a_i >0,\\
0 & \text{if } a_i=0,
\end{cases}\\
\beta_j \hook \left(x_{0}^{(a_0)} \dotsm x_{N_1}^{(a_{N_1})} \cdot y_{0}^{(b_0)} \dotsm y_{N_2}^{(b_{N_2})}\right) & =
\begin{cases}
x_{0}^{(a_0)} \dotsm x_{N_1}^{(a_{N_1})} \cdot y_{0}^{(b_0)} \dotsm y_i^{(b_i-1)} \dotsm y_{N_2}^{(b_{N_2})} & \text{if } b_i >0,\\
0 & \text{if } b_i=0,
\end{cases}
\end{align*}
and
$(\Theta\cdot\Psi)\hook F := \Theta\hook (\Psi\hook F)$
together with bilinearity
\[
(\Theta + \Psi)\hook (F+G) = \Theta\hook F + \Theta \hook G + \Psi\hook F + \Psi \hook G.
\]
\begin{defin}
For $F\in \Bbbk_{dp}[x_0,\ldots, x_{N_1}, y_0, \ldots, y_{N_2}]$ we define
$\Ann_{\PP^{N_1}\times \PP^{N_2}}(F):=
\set{\Theta \mid \Theta\hook F =0}$,
which is an ideal in  $\Bbbk[\alpha_0,\ldots,\alpha_{N_1},\beta_0,\ldots,\beta_{N_2}]$.
Define also the apolar algebra of $F$:
\[
\Apolar[F]:=\Bbbk[\alpha_0,\ldots,\alpha_{N_1},
\beta_0,\ldots,\beta_{N_2}]
/
\Ann_{\PP^{N_1}\times \PP^{N_2}}(F).
\]
\end{defin}

Note that for any $F$ as above and an invertible $\lambda \in \kk$ we have
$\Ann_{\PP^{N_1}\times \PP^{N_2}}(F) = \Ann_{\PP^{N_1}\times \PP^{N_2}}(\lambda F)$ and analogously
for $\Apolar[F]$.
In particular, instead of annihilators of elements of some vector space we can equally well talk about
the annihilators of the corresponding points in the projective space.

\begin{lemma}\label{lem_apolarity_on_PxP}
Let $p\in \PP(\Bbbk_{dp}[x_0,\ldots, x_{N_1}, y_0, \ldots, y_{N_2}]_{d,e})$ be a point and consider a  bihomogeneous ideal
$I \subset \Bbbk[\alpha_0,\ldots,\alpha_{N_1},\beta_0,\ldots,\beta_{N_2}]$.
Then
\begin{enumerate}
\item
\label{item_apolarity_for_ideals_in_PxP}
$I\subset \Ann_{\PP^{N_1}\times\PP^{N_2}}(p)$ if and only if $p \in \PP(I_{d,e}^{\perp})$,
\item
\label{item_apolarity_for_linear_spans_in_PxP}
In particular, if $I = J_Y$ is a saturated ideal
of a subscheme $Y \subset \PP^{N_1}\times\PP^{N_2}$,
then  $J_Y \subset \Ann_{\PP^{N_1}\times\PP^{N_2}}(p)$
if and only if $p$ is in the linear span of $Y$
reembedded via the linear system
$|\ccO_{\PP^{N_1}\times\PP^{N_2}}(d,e)|$,
\item
\label{item_symmetry_of_Hilbert_f}
$\dim\Apolar_{i,j} = \dim \Apolar_{d-i,e-j}$.
\end{enumerate}
\end{lemma}

The proof is no different than over $\CC$ \cite[Thm~1.4, Prop.~4.5]{galazka_multigraded_apolarity},
and also no different than single graded case, see \cite[Prop.~3.4]{nisiabu_jabu_cactus}.
This lemma has an interesting consequence for our setting.\begin{defin}
Suppose $X$ is a projective scheme, $(L_1, L_2)$ is a doubly standard graded.
For any
\[
p\in \PP(H^0(L_1^d\otimes L_2^e)^*) \subset \PP(\Bbbk_{dp}[x_0,\ldots, x_{N_1}, y_0, \ldots, y_{N_2}]_{d,e})
\]
define
$\Ann(p)=\Ann_{S_{L_1,L_2}}(p) \subset S_{L_1, L_2} = \Bbbk[\alpha_0,\ldots,\alpha_{N_1},\beta_0,\ldots,\beta_{N_2}]/I_X$
to be the ideal defined as $\Ann_{\PP^{N_1}\times \PP^{N_2}}(p)/I_X$.
\end{defin}
Note that in the setting of the definition we always have
$I_X \subset \Ann_{\PP^{N_1}\times \PP^{N_2}}(p)$ by Lemma~\ref{lem_apolarity_on_PxP}, since $X$ is
not contained in any hyperplane in $\PP(H^0(L_1^d\otimes L_2^e)^*)$.
In particular, we have
$
\Apolar  = S_{L_1,L_2}/ \Ann(p)
$
and the symmetry of Hilbert function (Lemma~\ref{lem_apolarity_on_PxP}\ref{item_symmetry_of_Hilbert_f})
applies equally well to the $S_{L_1, L_2}$ setting.
Further apolarity claims we need for $S_{L_1, L_2}$ are also analogous to Lemma~\ref{lem_apolarity_on_PxP}.

\begin{prop}[Apolarity lemma]\label{prop_apolarity_for_S_L1_L2}
  Suppose $X$ is a projective scheme, $(L_1, L_2)$ is a doubly standard graded, $d,e\geqslant 0$,
  with $(d,e) \neq (0,0)$, and let $p\in \PP(H^0(L_1^d\otimes L_2^e)^*)$.
Then:
\begin{enumerate}
\item \label{item_apolarity_for_ideals_in_S_L1_L2}
If $I\subset S_{L_1,L_2}$ is any bihomogeneous ideal,
then $I \subset \Ann(p)$ if and only if $p \in \PP(I_{d,e}^{\perp})$.
\item \label{item_apolarity_for_linear_spans_in_S_L1_L2}
If $R\subset X$ is a closed subscheme of $X$,
then $J_R \subset\Ann(p)$ if and only if $p$ is in the linear span of $R$ reembedded
into $X\subset \PP(H^0(L_1^d\otimes L_2^e)^*)$.
\end{enumerate}
\end{prop}
\begin{prf}
Let
$\hat I \subset \Bbbk[\alpha_0,\ldots,\alpha_{N_1},
\beta_0,\ldots,\beta_{N_2}]$
be the preimage of $I$ under the quotient map.
Then the statement of \ref{item_apolarity_for_ideals_in_S_L1_L2} follows directly
from $I_X \subset \hat I$, $I_X \subset \Ann_{\PP^{N_1}\times \PP^{N_2}}(p)$ and
Lemma~\ref{lem_apolarity_on_PxP}\ref{item_apolarity_for_ideals_in_PxP}.

Similarly, define $\hat{J}_Y$ to be the preimage under the same quotient map of $J_Y$.
Then $\hat{J}_Y$ is the saturated bihomogeneous ideal of $Y$ in $\PP^{N_1}\times \PP^{N_2}$ and
the claim of \ref{item_apolarity_for_linear_spans_in_S_L1_L2}
follows from Lemma~\ref{lem_apolarity_on_PxP}\ref{item_apolarity_for_linear_spans_in_PxP}.
\end{prf}

\begin{lemma}\label{lem_non_zero_apolar_algebra}
  Suppose $X$ is a projective scheme, $(L_1, L_2)$ is a doubly standard graded, $d,e\geqslant 0$,
  with $(d,e) \neq (0,0)$, and let $p\in \PP(H^0(L_1^d\otimes L_2^e)^*)$.
Then for any integers $i, j$,
the $(i,j)$ grading of the apolar algebra $\Apolar_{i,j} \neq 0$ if and only if $0\leqslant i \leqslant d$
and $0\leqslant j \leqslant e$.
\end{lemma}

\begin{prf}
  Since $\Apolar_{i,j} = (S_{L_1,L_2}/\Ann(p))_{i,j}$ and the negative gradings of the ring $S_{L_1,L_2}$ are all zero,
  which proves the claim for $i< 0$ and for $j< 0$.
  By the symmetry in Lemma~\ref{lem_apolarity_on_PxP}\ref{item_symmetry_of_Hilbert_f} this also proves
  the claim for $i >d$ and for $j> e$.

  If $i=d$ and $j=e$, then not all sections from $S_{L_1,L_2}$ annihilate $p$, as otherwise $p$ would be $0$, which
  is not a point of the projective space.
Thus $\Apolar_{d,e}$ and also $\Apolar_{0,0}$ (by symmetry again) are non-zero.
On the other hand, $\Apolar_{1,0}$ and $\Apolar_{0,1}$ generate the algebra $\Apolar$, thus by downward induction
we prove that $\Apolar_{i,j} \neq 0$ for any of the remaining $i,j$:
Indeed, if by contradiction $\Apolar_{i,j} = 0$, then also $\Apolar_{i+1,j} = \Apolar_{i,j+1} = 0$ by
Lemma~\ref{lem_doubly_standard_graded_algebra}, which eventually is a contradiction with $\Apolar_{d,e}\neq 0$.
\end{prf}

\section{Cactus varieties}
\label{sec_cacti}

In this section we prove our main result, Theorem~\ref{thm_main_intro_simpler}.
For this purpose, throughout this section  we work with a fixed projective scheme $X$
and a positive integer $r$.

Suppose $L$ is a very ample line bundle on $X$, and $A$ and $B$ are two line bundles on $X$ such that
$L=A\otimes B$. In this situation we will say that this expression $L=A\otimes B$ is a \emph{splitting} of $L$.
The multiplication of sections
$H^0(X,A) \otimes H^0(X,B)  \to H^0(X,L)$ is a linear map that can be represented by the following matrix
$M=M_{A,B}$ with entries in $H^0(X,L)$.
Pick bases $\fromto{a_{1}}{a_{\dim H^0(X,A)}}$ and $\fromto{b_{1}}{b_{\dim H^0(X,B)}}$ of the spaces
of sections of $H^0(X,A)$ and $H^0(X,B)$ respectively and let the $(i,j)$-th entry of $M$ is $a_i\cdot b_j\in H^0(X,L)$.
This matrix depends on the choices of the bases, but different choices lead to
$\GL(H^0(X,A)) \times \GL(H^0(X,B))$-equivalent matrices.
In particular, for any integer $r'$ the ideal generated by $(r'+1)\times (r'+1)$ minors of $M$ does not depend on the choices.
Thus, simplifying, we will not mention  the choices of bases in our further considerations $M$,
we will simply say that $M$ is the \emph{matrix associated to the splitting}  $L=A\otimes B$.

For the induction purposes, we will need to deal with somewhat trivial (but deserving different treatment than general) case
$r=0$ of Theorem~\ref{thm_main_intro_simpler}, which we deal with here.

\begin{lemma}
\label{lem_case_r_equal_to_0}
With the notation as above,
suppose $(A,B)$ is double standard graded.
Then the $1\times 1$-minors
of $M$ generate the maximal homogeneous
ideal of $S_L$.
\end{lemma}
\begin{prf}
The $1\times 1$-minors of a matrix are just its entries.
Since the multiplication of sections  $H^0(X,A) \otimes H^0(X,B)  \to H^0(X,L)$
is surjective by the definition of double standard graded, the entries of $M$, which are $a_i\cdot b_j$,
generate $H^0(X,L)$.   Since $L=A\otimes B$, thus $L$ is also standard graded, so that the ideal
$(H^0(X,L))\subset S_L$ contains all $(S_L)_d$ for any $d>0$, and thus $S_L/(H^0(X,L)) \simeq (S_L)_0 = \kk$,
proving the claim.
\end{prf}

Now we want to compare the apolarity theory from \S\ref{sec_apolarity} to multiplication of sections,
which again is standard knowledge in the case of apolarity in a polynomial ring.
Suppose that $X$ is a projective scheme and $(L_1,L_2)$ is doubly standard graded.
First, we observe that for any $i,j$ the standard duality map
$H^0(L_1^i\otimes L_2^j) \otimes H^0(L_1^i\otimes L_2^j)^* \to \kk $ coincides with $\hook$ in this degree,
$\Theta \otimes F \mapsto \Theta \hook F\in \kk=(S_{L_1,L_2})_{0,0}$.
Then for fixed integers $d,e,i,j$ set $A := L_1^{i} \otimes L_2^{j}$ and $B:=L_1^{d-i} \otimes L_2^{e-j}$.
Pick a point $p\in \PP(H^0(L_1^d\otimes L_2^e)^*)$,
and consider the evaluation $M_{A,B}(p)$,
for a matrix $M_{A,B}$ as above.
For fixed bases, the evaluation is well defined up to a rescaling, depending on the lifting of $p$ from
the projective space to the affine space.
But the rank of this matrix $M_{A,B}(p)$ is independent of the rescaling.

\begin{lemma}
\label{lem_dim_of_apolar_is_rk}
Suppose $X$ is a projective scheme, $(L_1, L_2)$
is a doubly standard graded,
$d,e\geqslant 0$, with $(d,e) \neq (0,0)$,
and let $p\in \PP(H^0(L_1^d\otimes L_2^e)^*)$.
Then for two integers $i, j$ such that $0\leqslant i \leqslant d$ and $0 \leqslant j \leqslant e$ and $(i,j)\neq (0,0), (d,e)$,
we have $\dim_{\kk} ( \Apolar_{i,j}) =
\rk M_{A,B}(p)$, where $A := L_1^{i} \otimes L_2^{j}$ and $B:=L_1^{d-i} \otimes L_2^{e-j}$ as above.
\end{lemma}

\begin{prf}
The assumptions on $d,e, i,j$ assure that
$(S_{L_1,L_2})_{i,j} = H^0(A)$ and
$(S_{L_1,L_2})_{d-i,e-j} = H^0(B)$,
and analogously for $A\otimes B$
and $(d,e)$-th grading.
The matrix $M_{A,B}(p)$ arises from composing the multiplication of sections and evaluation at $p$:
\[
(S_{L_1,L_2})_{i,j} \otimes (S_{L_1,L_2})_{d-i,e-j} \to (S_{L_1,L_2})_{d,e} \stackrel{p}{\to} \kk.
\]
Thus rank of $M_{A,B}(p)$ is equal to rank of the following linear map:
\begin{align*}
(S_{L_1,L_2})_{i,j}& \to  ((S_{L_1,L_2})_{d-i,e-j})^*,\\
\Theta& \mapsto (\Psi \mapsto (\Psi\cdot \Theta) (p) ).
\end{align*}
By the comparison between apolarity and duality, we have:
$(\Psi\cdot \Theta) (p) = (\Psi\cdot \Theta) \hook p= \Psi \hook(\Theta \hook p)$.
Thus explicitly, the linear map above takes $\Theta$ to
$\Theta \hook p\in ((S_{L_1,L_2})_{d-i,e-j})^*$.
Then $\Theta$ is in the kernel of this linear map if and only if $\Theta\in \Ann(p)$, and
\[
\rk M_{A,B}(p) = \dim_{\kk} (S_{L_1,L_2})_{i,j}
- \dim_{\kk} \Ann(p)_{i,j} = \dim_{\kk} \Apolar_{i,j}.
\]
\end{prf}

Recall, that for any embedding $X\subset \PP^N$, the $r$-th cactus variety of $X$ is defined as:
\[
\cactus{r}{X}=
\overline{\bigcup
\set{\langle R \rangle  \mid R \subset X,
\text{ $R$ is a subscheme of degree
at most $r$}}} \subset \PP^N.
\]

We split the proof of Theorem~\ref{thm_main_intro_simpler} into two parts.
In the first part, we list conditions on line bundles $A$ and $B$ that guarantee that the cactus variety of $X$
embedded by $L= A\otimes B$ is defined by minors of the matrix as above.
This part is algebraic and mimics the proofs in  \cite[\S5]{nisiabu_jabu_cactus}, but adapted
to this more general situation. In the second part, we show that for a sufficiently ample line bundle $L$
there exists a splitting $L=A\otimes B$ satisfying these conditions.
This proof is largely combinatorial, and in the end it boils down to showing that there exists a lattice point
in some region of $\RR^n$ bounded by a (possibly infinite) collection of linear inequalities.

\subsection{Algebraic part: exploiting apolarity in double section ring}
\begin{interruptedenumerate}

\begin{thm}
\label{thm_main_with_splitting}
Pick integers $k$ and $d$  and line bundles $L_1$ and $L_2$ on $X$ such that 
\begin{insideenumerate}
\item
\label{item_thm_main_assumption_on_r_k_d}
$r \leqslant  k \leqslant d-r$,
\item
\label{item_thm_main_assumption_on_L1_generated}
$L_1$ is $(k+1)$-embedding in the sense of Definition~\ref{def_d_embedding},
\item
\label{item_thm_main_assumption_on_L1_L2_doubly_stndrd}
the pair $L_1, L_2$ is doubly standard graded in the sense of Definition~\ref{def_doubly_standard_graded}, and
\item
\label{item_thm_main_assumption_on_the_other_pair_doubly_stndrd}
the multiplication of sections
$H^0\left(X,L_1^{\otimes (d-r)} \otimes L_2^*\right)\otimes H^0(X,L_2) \to H^0\left(X,L_1^{\otimes (d-r)}\right)$ is surjective.
\end{insideenumerate}
Set $L= L_1^{\otimes d} \otimes L_2$ and  consider the embedding of $X$ via the complete linear system $H^0(X, L)$.
Then the cactus variety $\cactus{r}{X}$ is set-theoretically determined by $(r+1)\times (r+1)$-minors
of the matrix with linear entries defined by the splitting
$L=A\otimes B$ with $A= L_1^{\otimes k}$ and $B=L_1^{\otimes d-k} \otimes L_2$.
\end{thm}

For the proof, we need another lemma.
As a consequence of the first one we will be able to use downward induction on $r$ in the proof of the theorem.
If the assumptions
\ref{item_thm_main_assumption_on_r_k_d}--\ref{item_thm_main_assumption_on_the_other_pair_doubly_stndrd}
are satisfied for some value of $r$, then they are also satisfied for any smaller, non-negative value of $r$.

\begin{lemma}
\label{lem_higher_power_of_L1_is_fine}
In the notation of Theorem~\ref{thm_main_with_splitting}
suppose assumptions~\ref{item_thm_main_assumption_on_L1_L2_doubly_stndrd} and \ref{item_thm_main_assumption_on_the_other_pair_doubly_stndrd} hold.
Then also:
\begin{insideenumerate}
\item
\label{item_thm_main_assumption_on_the_other_pair_for_higher_d_doubly_stndrd}
the multiplication of sections $H^0\left(X,L_1^{\otimes e} \otimes L_2^*\right)
\otimes
H^0(X,L_2)
\to
H^0\left(X,L_1^{\otimes e}\right)$ is surjective for any $e\geqslant d-r$.
\end{insideenumerate}
\end{lemma}
\end{interruptedenumerate}

\begin{prf}
If $e= d-r$ then there is nothing to prove.
By induction, it is enough show the claim for $e=d-r+ 1$.
Below, for brevity we skip $X$ in the notation for $H^0(\dots)$.
The map of interest
\[
H^0\left(L_1^{\otimes (d+1-r)} \otimes L_2^*\right)
\otimes
H^0(L_2)
\to
H^0\left(L_1^{\otimes (d+1-r)}\right)
\]
fits as the bottom map into a commutative diagram:
\[
\begin{tikzcd}
H^0(L_1)\otimes H^0\left(L_1^{\otimes (d-r)} \otimes L_2^*\right)
\otimes
H^0(L_2)\arrow[d] \arrow[r] &      H^0(L_1)\otimes H^0\left(L_1^{\otimes (d-r)}\right)\arrow[d]\\
H^0\left(L_1^{\otimes (d+1-r)} \otimes L_2^*\right)
\otimes
H^0(L_2) \arrow[r]      & H^0\left(L_1^{\otimes (d+1-r)}\right).
\end{tikzcd}
\]
The top map is surjective by
\ref{item_thm_main_assumption_on_the_other_pair_doubly_stndrd}, and
the right map is surjective by
\ref{item_thm_main_assumption_on_L1_L2_doubly_stndrd}.
Therefore also the bottom map must be surjective, as claimed.
\end{prf}

\begin{prf}[ of Theorem~\ref{thm_main_with_splitting}]
Let $M$ be the matrix of linear forms on   $\PP(H^0(X, L)^*)$
arising from the splitting $L=A\otimes B$.
The inclusion   $\cactus{r}{X}\subset     \reduced{\set{\rk M\leqslant r}}$
follows from \cite[Thm~5]{galazka_vb_cactus}    or \cite[Thms~1.17, 1.18]{galazka_phd}.
Thus it remains to show the opposite inclusion.
That is, pick any point $p \in \reduced{\set{\rk M\leqslant r}}\subset \PP(H^0(X, L)^*)$ and
we aim to construct a subscheme $R \subset X$ of degree at most $r$ such that $p\in \linspan{R}$.
More precisely,  we will construct a scheme $R$ of degree equal to $\rk M(p)$.
Without loss of generality   (potentially replacing $r$ by a lower value)
we may assume $r = \rk M(p)$.
Note that, decreasing $r$ does not violate any of the
assumptions~\ref{item_thm_main_assumption_on_r_k_d}--\ref{item_thm_main_assumption_on_L1_L2_doubly_stndrd},
and \ref{item_thm_main_assumption_on_the_other_pair_doubly_stndrd} still holds thanks to
Lemma~\ref{lem_higher_power_of_L1_is_fine} and
its item~\ref{item_thm_main_assumption_on_the_other_pair_for_higher_d_doubly_stndrd}.

We proceed to define several bihomogeneous ideals in $S_{L_1, L_2}$.
Consider $\Ann(p)\subset S_{L_1, L_2}$, which is a cofinite and bihomogeneous ideal.
Recall the apolar algebra $\Apolar = S_{L_1, L_2}/\Ann(p)$, which is a bigraded finite $\kk$-algebra
with a symmetry $\dim_{\kk}\Apolar_{(i, j)} = \dim_{\kk}\Apolar_{(d-i, 1-j)}^*$ by
Lemma~\ref{lem_apolarity_on_PxP}\ref{item_symmetry_of_Hilbert_f}.
Moreover, $\Apolar_{(i, j)} \neq 0$ if and only if $0\leqslant i \leqslant d$ and $0\leqslant j \leqslant 1$ by
Lemma~\ref{lem_non_zero_apolar_algebra}.
By our assumptions and Lemma~\ref{lem_dim_of_apolar_is_rk}:
\[
\dim  \Apolar_{(k, 0)}= r = \rk M(p).
\]
By Macaulay bound, since $k\geqslant r$, the Hilbert function on $(\ast, 0)$-line must be non-increasing
from $(k,0)$ onwards, that is:
\[
\dim  \Apolar_{(k', 0)} \geqslant \dim  \Apolar_{(k'+1, 0)} \text{ for all } k'\geqslant k.
\]
But $\dim  \Apolar_{(d, 0)} \geqslant 1$ (so the drop is at most $r-1$) and
there are $d-k \ge r$ steps between $k$ and $d$.
So there must exist $k_0$ with $k \leqslant k_0 \leqslant d-1$ such that
$\dim  \Apolar_{(k_0, 0)} = \dim  \Apolar_{(k_0+1, 0)}\geqslant 1$.
Choose $k_0$ to be the minimal integer satisfying all of the above conditions,
and let $r_0 = \dim  \Apolar_{(k_0, 0)}\leqslant r$.
Eventually, we will prove that $r_0=r$ and thus $k_0=k$.

Now define the ideal    $I\subset S_{L_1}$   generated by the first $k_0$ gradings of
$ \Ann(p)_{(\ast, 0)}$:
\[
I := \left(\Ann(p) \cap S_{L_1}\right)_{\leqslant k_0} =
\left( \textstyle \bigoplus_{i=0}^{k_0}  \Ann(p)_{(i,0)} \right).
\]
Express   $S_{L_1} = \kk[\fromto{\alpha_0}{\alpha_{N_1}}] / I_X$
where  $I_X \subset \kk[\fromto{\alpha_0}{\alpha_{N_1}}]$ is the  homogeneous ideal defining
$X\subset \PP^{N_1} =\PP(H^0(X,L_1)^*)$.
By the assumption of the theorem $I_X$ is generated in degrees at most $k+1$.

Let $\widehat{I}\subset \kk[\fromto{\alpha_0}{\alpha_{N_1}}]\left/ I_X\right.$ be the lift of $I$ to
the polynomial ring $\kk[\fromto{\alpha_0}{\alpha_{N_1}}]$.
Let $\widehat{I}^{\sat}$ be the saturation of $\widehat{I}$ with respect to the maximal homogeneous
ideal $(\fromto{\alpha_0}{\alpha_{N_1}})$.
Note that $I_X\subset \widehat{I} \subset \widehat{I}^{\sat}$
Also let $I^{\sat} \subset S_{L_1}$ be the descended ideal, so that $\widehat{I^{\sat}} = \widehat{I}^{\sat}$.
Thus $\widehat{I}$ is generated by the generators of $I_X$ and by lifts of the generators of $I$.
In particular, $\widehat{I}$ is generated in degrees at most $k_0+1$.
However, the growth of $S_{L_1}/I = \kk[\fromto{\alpha_0}{\alpha_{N_1}}]\left/ \widehat{I}\right.$
is maximal possible from degree $k_0$ to $k_0+1$, thus there cannot be any minimal generator of
$\widehat{I}$ in degree $k_0+1$.
Thus $\widehat{I}$ is generated in degrees at most $k_0$.
Therefore, by Lemma~\ref{lem_Gotzmann_in_single_graded}
both algebras $\kk[\fromto{\alpha_0}{\alpha_{N_1}}]\left/\widehat{I}\right.$ and
$ \kk[\fromto{\alpha_0}{\alpha_{N_1}}]\left/\widehat{I}^{\sat}\right.$ have their Hilbert functions constant
and equal to $r_0$ from $k_0$ onwards and therefore $(\widehat{I}^{\sat})_{\geqslant k_0} = \widehat{I}_{\geqslant k_0}$.
In particular, $\widehat{I}^{\sat}$ is generated in degrees at most $k_0$.
Note that
$
\dim\left(\kk[\fromto{\alpha_0}{\alpha_{N_1}}]\left/\widehat{I}^{\sat}\right)_{i}\right. \leqslant r_0$ for any $i$.
Furthermore, $ I^{\sat}\subset \Ann(p)\cap S_{L_1}$. Thus, looking at degree $k\leqslant k_0$ part, we must
have $r_0\geqslant r$, thus $r_0=r$ and $k_0=k$.
Define $R=Z(I^{\sat})$ to be the subscheme of $X$ of degree $r$.

Now we switch back our attention to ideals in the double section ring $S_{L_1,L_2}$.
Let $J:=J_R\subset S_{L_1,L_2}$ be the bihomogenous saturated (with respect to both ideals
$\left((S_{L_1,L_2})_{(1,0)}\right)$ and
$\left((S_{L_1,L_2})_{(0,1)}\right)$)
ideal of $R\subset X$.
Let $\widehat{J} \subset
\kk[\fromto{\alpha_0}{\alpha_{N_1}}, \fromto{\beta_0}{\beta_{N_2}}]$ be the lift of $J$ to the bigraded polynomial ring.
This ideal is also saturated, and corresponds to the ideal sheaf $\ccJ\subset \ccO_{\PP^{N_1}\times \PP^{N_2}}$
defining $R\subset \PP^{N_1}\times \PP^{N_2}$.
The Hilbert function $h_J\colon \ZZ^2 \to \NN$
of
$\kk[\fromto{\alpha_0}{\alpha_{N_1}},
\fromto{\beta_0}{\beta_{N_2}}]\left/\widehat{J}\right.$
(which is equal to $S_{L_1,L_2}/J$)
is bounded from above by the degree of $R$,
that is $r$.
It is also non-decreasing, $h_J(i,j) \leqslant h_J(i+1,j)$ and $h_J(i,j) \leqslant h_J(i,j+1)$ by
Lemma~\ref{lem_Gotzmann_and_Hilbert_function}\ref{item_monotone_nonzerodivisor_exists}.
We also have that $h_J(r-1,0) =r$
as $J$ agrees with $I$ in the $(\ast,0)$-gradings.
Thus $\widehat{J}$ is $(r,1)$-regular in the sense of \cite[Def.~4.1]{maclagan_smith_multigraded_regularity}.
By \cite[Thm.~5.4]{maclagan_smith_multigraded_regularity}
there are no minimal generators of $\widehat{J}$ in degrees $(r',1)$ for $r'>r$,
and by previous considerations on $I$ we know that there are no minimal generators in degrees $(r',0)$ for $r'>r$.

To show that $p\in \linspan{R}$ we must show that $J_{(d,1)} \subset \Ann(p)_{(d,1)}$ and apply
Proposition~\ref{prop_apolarity_for_S_L1_L2}\ref{item_apolarity_for_linear_spans_in_S_L1_L2}.
By the above considerations,
\[
J_{(d,1)} = J_{(r,1)} \cdot (S_{L_1, L_2})_{(d-r,0)} =  J_{(r,1)} \cdot H^0(X,L_1^{d-r}).
\]
By our assumption~\ref{item_thm_main_assumption_on_the_other_pair_doubly_stndrd}
we have  $H^0(X,L_1^{d-r}) = H^0(X,L_1^{d-r}\otimes L_2^*) \cdot H^0(X,L_2)$.
Thus we have:
\begin{align*}
J_{(d,1)} &= J_{(r,1)} \cdot H^0(X,L_1^{d-r}) = J_{(r,1)} \cdot H^0(X,L_1^{d-r}\otimes L_2^*) \cdot H^0(X,L_2)  \\
& \subset J_{(d,0)} \cdot H^0(X,L_2) = I_{d} \cdot H^0(X,L_2) \\
& \subset \Ann(p)_{(d,0)} \cdot H^0(X,L_2) \subset \Ann(p)_{(d,1)},
\end{align*}
as claimed.
Thus indeed $p\in \linspan {R}$
and $p\in \cactus{r}{X}$.
\end{prf}

\subsection{Combinatorial part: lattice points and inequalities}

The goal of this section is to show:

\begin{thm}
\label{thm_main_combinatorial_part}
Suppose $L$ is a sufficiently ample line bundle on $X$.
Then there exist integers $d$, $k$, and line bundles $L_1$, $L_2$ satisfying
properties~\ref{item_thm_main_assumption_on_r_k_d}--\ref{item_thm_main_assumption_on_the_other_pair_doubly_stndrd}
of Theorem~\ref{thm_main_with_splitting} and such that $L = L_1^{\otimes d} \otimes L_2$.
\end{thm}

In the first step we reduce the claim about existence of some line bundles (thus working in $\Pic(X)$) to
an analogous claim in the lattice $N^1(X)\simeq \ZZ^{\rho}$.
For this purpose we denote by $\QDSG\subset \Pic(X)$ (for Quadratically generated and Doubly Standard Graded)
a set of line bundles such that
\begin{itemize}
\item $\QDSG = \set{L_0\otimes \Lambda \mid \Lambda \text{ is nef line bundle}}$ for some line bundle $L_0$,
\item for any $L_1, L_2 \in \QDSG$ the pair $(L_1, L_2)$ is doubly standard graded.
\item for any $L_1 \in \QDSG$ the line bundle  $L_1$ is $2$-generated.
\end{itemize}
Such set $\QDSG$ exists by Corollary~\ref{cor_standard_graded_for_suff_ample}\ref{item_doubly_standard_graded}
and by \cite[Thm~1.1]{sidman_smith_linear_det_eq_for_all_proj_schemes}
(combined with Lemma~\ref{lem_finitely_many_properties_of_sufficiently_ample_lbs} and
Proposition~\ref{prop_nef_definition_of_sufficiently_ample}).
Note that $\QDSG$ is not uniquely defined, but any of them will do the job.

By $D_0\in N^1(X)$ denote the numerical class of $L_0$.
Thus for any line bundle $L_1$,
we have $L_1\in \QDSG$ if and only if the numerical class of $L_1$ is in $D_0+\Nef(X)$.

\begin{prop}
\label{prop_reduction_to_lattice_conditions}
Let $L$ be an ample line bundle and let $D\in N^1(X)$ be its numerical class.
Suppose there exist integer $d$ and lattice points $D_1$, $D_2$ such that:
\renewcommand{\theenumi}{\textnormal{\textbf{\textsc{(\alph{enumi})}}}}
\begin{enumerate}
\item
\label{item_prop_on_reduction_to_lattice_assumption_on_r_k_d}
$d\geqslant 2r$,
\item
\label{item_prop_on_reduction_to_lattice_equality_D_D1_D2}
$D=d D_1 + D_2$,
\item
\label{item_prop_on_reduction_to_lattice_divisors_in_QDSG}
all of $D_1$, $D_2$ and
$(d-r) D_1 - D_2$
are contained in $D_0+\Nef(X)$.
\end{enumerate}
\renewcommand{\theenumi}{(\roman{enumi})}
Then there exist integer $k$ and line bundles $L_1$, $L_2$ that (together with $d$ above) satisfy properties
\ref{item_thm_main_assumption_on_r_k_d}--\ref{item_thm_main_assumption_on_the_other_pair_doubly_stndrd}
of Theorem~\ref{thm_main_with_splitting} and such that $L = L_1^{\otimes d} \otimes L_2$.
\end{prop}
\begin{prf}
To satisfy \ref{item_thm_main_assumption_on_r_k_d}
pick any $k$ such that $r\leqslant k \leqslant d-r$,
which is possible by \ref{item_prop_on_reduction_to_lattice_assumption_on_r_k_d}.
Let $L_1$ be any line bundle whose numerical class is $D_1$ and define $L_2 = L\otimes (L_1^*)^{\otimes d}$.
Then the numerical class of $L_2$ is $D_2$
by \ref{item_prop_on_reduction_to_lattice_equality_D_D1_D2}.
Moreover, by \ref{item_prop_on_reduction_to_lattice_divisors_in_QDSG} we have:
\begin{itemize}
\item $L_1, L_2 \in \QDSG$ so that \ref{item_thm_main_assumption_on_L1_generated}
and  \ref{item_thm_main_assumption_on_L1_L2_doubly_stndrd} hold,
\item $L_1^{\otimes (d-r)} \otimes L_2^* \in \QDSG$ so
that \ref{item_thm_main_assumption_on_the_other_pair_doubly_stndrd} holds,
\end{itemize}
concluding the proof.
\end{prf}

Now we work in the lattice $N^1(X)$ and the corresponding vector space $N^1_{\RR}(X)$.
Assuming \ref{item_prop_on_reduction_to_lattice_assumption_on_r_k_d} and
\ref{item_prop_on_reduction_to_lattice_equality_D_D1_D2} we translate
the conditions \ref{item_prop_on_reduction_to_lattice_divisors_in_QDSG} of the proposition into some
translated cone conditions.
We may assume that $D$ is a class of a sufficiently ample line bundle, which we will write
$D = \lambda D_0+ \Delta$ for $\Delta\in \Nef(X)$ and a sufficiently large integer $\lambda$ (to be determined later).
Write also $D_1= D_0+ \Delta_1$.
Then
\begin{align*}
D_2 & = D - d D_1 = (\lambda-d) D_0 + \Delta - d \Delta_1,\\
(d-r) D_1 - D_2 & = (d-r)(D_0+ \Delta_1) - (\lambda-d) D_0 - \Delta + d \Delta_1\\
&= (2d - r-  \lambda) D_0 - \Delta - (2d-r) \Delta_1.
\end{align*}
Therefore, in this notation \ref{item_prop_on_reduction_to_lattice_divisors_in_QDSG} is equivalent to
the following three classes being in $\Nef(X)$:
\begin{equation}
\label{equ_nef_inequalities_for_Deltas_and_Ds}
\Delta_1, \quad
(\lambda-d-1) D_0 + \Delta - d \Delta_1, \quad
(2d - r-  \lambda-1) D_0 - \Delta + (2d-r) \Delta_1.
\end{equation}

In order to proceed, we choose  an ample class $H\in \Amp(X)$ so that
$\Nef(X)$ and $\Amp(X)$ are defined by a compact set of inequalities $S = \overline{NE(X)}\cap (H=1) \subset N^1(X)^*$:
\[
\Nef(X) = \set{E \mid \forall_{\phi\in S}  \ \phi(E)\geqslant 0 } \text{ and }
\Amp(X) = \set{E \mid \forall_{\phi\in S}  \ \phi(E)> 0}.
\]
As in Lemma~\ref{lem_minD_and_maxD}, for each $E\in N^1(X)$ we have a well defined
\[
\minD[E] =  \min\set{\phi(E) \mid \phi \in S}.
\]

We also choose any $\ZZ$-basis $\fromto{e_1}{e_{\rho}}$ of $N^1(X)$, and consider the cube $C$, whose
vertices are $\pm e_1 \pm \dots \pm e_{\rho}$.

\begin{lemma}
\label{lem_exists_integral_point_in_Nef}
Suppose $D_0 \in \Amp(X)$ and $\lambda$ is an integer satisfying the following conditions:
\begin{itemize}
\item $\lambda \geqslant 2d-r-1$, and
\item $\lambda \geqslant \dfrac{3d-r}{d-r} - \dfrac{d(2d-r)\minD[v]}{(d-r)\minD[D_0]}$ for any choice
  of the vertex $v=\pm e_1 \pm \dots \pm e_{\rho}$ of the cube $C$.
\end{itemize}
Then for any $\Delta\in\Nef(X)$
there exists an integral lattice point
$\Delta_1\in N^1(X)$ such that all three classes from \eqref{equ_nef_inequalities_for_Deltas_and_Ds} are in $\Nef(X)$.
\end{lemma}

\begin{rem}
\label{rem_sublattice_is_good_enough}
If we force appropriately stronger bound on $\lambda$
in Lemma~\ref{lem_exists_integral_point_in_Nef}
we can also insist that $\Delta_1$
is in a finer sublattice of $N^1_{\RR}(X)$.
In fact, even if $N'\subset N^1(X)$ is any sublattice (not necessarily of full dimension) such that
$N' \cap \Amp(X)\ne \emptyset$ and both $D_0$ and $\Delta$ are in $N'_{\RR}$, then there exists $\Delta_1\in N'$
such that all classes from \eqref{equ_nef_inequalities_for_Deltas_and_Ds} are in $\Nef(X)\cap N'_{\RR}$.
The only modification in the proof (and the statement) is that we must replace the basis of $N^1(X)$ with a basis of $N'$ (and correspondingly, the vertices of the cube, and the bound on $\lambda$).
\end{rem}

\begin{prf}
Using elementary linear operations we translate the conditions on nefness of three classes in
\eqref{equ_nef_inequalities_for_Deltas_and_Ds} into three conditions on $\Delta_1$:
\begin{align*}
\Delta_1 & \in \Nef(X), \\
\Delta_1 & \in \left(\frac{\lambda-1}{d} -1\right) D_0 + \frac{1}{d}\Delta
- \Nef(X), \text{ and}\\
\Delta_1 &\in \left(\frac{\lambda +1}{2d-r} -1\right) D_0 + \frac{1}{2d-r}\Delta + \Nef(X).
\end{align*}
Since $\lambda\geqslant 2d-r-1$, and $d\geqslant 2r$ both
coefficients  $\frac{\lambda +1}{2d-r} -1$ and $\frac{1}{2d-r}$  of the third class are non-negative, and since
both $D_0$ and $\Delta$ are in $\Nef(X)$, the first condition follows from the third (so the first one is redundant).
Therefore, we must show that there is a lattice point $\Delta_1$ in the region bounded by a translated $\Nef(X)$
and by a translated $-\Nef(X)$. To obtain this goal,
we use Lemma~\ref{lem_condition_on_lattice_point_in_cones}.
We have to show that for each vertex
$v=\pm e_1 \pm \dots \pm e_{\rho}$ of $C$ the following class is nef
\begin{align*}
&
\left(\left(\frac{\lambda-1}{d} -1\right) D_0 + \frac{1}{d}\Delta\right) - \left(\left(\frac{\lambda +1}{2d-r} -1\right) D_0
+ \frac{1}{2d-r}\Delta\right) + v \\
= \ & \left(\frac{\lambda-1}{d} -\frac{\lambda +1}{2d-r}  \right) D_0 + \left( \frac{1}{d}- \frac{1}{2d-r}\right)\Delta + v.
\end{align*}
Since $d\geqslant 2r$, we must have
$\frac{1}{d} -\frac{1}{2d-r} \geqslant 0$.
Thus $\left(\frac{1}{d} -\frac{1}{2d-r}\right)\Delta$ is always nef for any choice of nef $\Delta$.
Therefore to conclude the proof it is enough to  show that the following class is nef:
\begin{align*}
\left(\frac{\lambda-1}{d} -\frac{\lambda +1}{2d-r} \right) D_0 + v
& =  \left(\frac{ (2d-r)(\lambda-1)- d(\lambda +1) }{d(2d - r)}  \right) D_0 + v\\
&=
\left(\frac{ (d-r)\lambda  - 3d+r }{d(2d - r)}  \right) D_0 + v.
\end{align*}
To check nefness we apply \eqref{equ_closed_cone_by_its_inequalities} and \eqref{equ_open_cone_by_its_inequalities}
for $\sigma = \Nef(X)$. Take any $\phi\in S$ and apply $\phi$ to the above class:
\begin{align*}
\phi\left(\left(\frac{ (d-r)\lambda  - 3d+r }{d(2d - r)}  \right) D_0 + v\right) &=   \left(\frac{ (d-r)\lambda  - 3d+r }{d(2d - r)}  \right) \phi(D_0) + \phi(v)\\
&\geqslant  \left(\frac{ (d-r)\lambda  - 3d+r }{d(2d - r)}  \right) \minD[D_0] + \minD[v] \geqslant 0.
\end{align*}
The last inequality ($\dots \geqslant 0$) follows from $\minD[D_0] > 0$, and our assumption on $\lambda$.
This concludes the proof  of the lemma.
\end{prf}

\begin{prf}[ of Theorem~\ref{thm_main_combinatorial_part}]
Pick any integer $d\geqslant 2r$
and any integer $k$, such that $r\leqslant k\leqslant d-r$.
Let $L_0$ be the very ample line bundle in the definition of $\QDSG  = L_0 \otimes (\text{nef line bundles})$
and let $D_0$ be the numerical class
of $L_0$.
Finally, let $\lambda$ be an integer satisfying the two inequalities in Lemma~\ref{lem_exists_integral_point_in_Nef}.

The sufficiently ample condition for $L$
that we are going to use is
$L = L_0^{\otimes \lambda}\otimes \Lambda$
with $\Lambda$ nef.
Let $\Delta\in \Nef(X)$ be the numerical class of $\Lambda$.
By Lemma~\ref{lem_exists_integral_point_in_Nef}
there exists an integral class $\Delta_1\in N^1(X)$
such that all numerical classes in~\eqref{equ_nef_inequalities_for_Deltas_and_Ds} are in $\Nef(X)$.
Equivalently, there exist lattice points $D_1$ and $D_2$ in $N^1(X)$ satisfying assumptions
of Proposition~\ref{prop_reduction_to_lattice_conditions}.
Thus this proposition implies that there exist line bundles $L_1$, $L_2$ such that assumptions of
Theorem~\ref{thm_main_with_splitting} are satisfied, concluding this proof.
\end{prf}

\begin{prf}[ of Theorem~\ref{thm_main_intro_simpler}]
First assume $r=r'$. Let $L$ be a sufficiently ample line bundle on $X$.
By Theorem~\ref{thm_main_combinatorial_part}
there exist integers $d$, $k$, and line bundles $L_1$ and $L_2$ such that $L=L_1^{\otimes d} \otimes L_2$
and they satisfy the assumptions of Theorem~\ref{thm_main_with_splitting}.
Thus by this theorem, $\cactus{r}{X}$ is defined by the desired minors, as claimed

Note that taking into account Lemma~\ref{lem_higher_power_of_L1_is_fine} conditions
\ref{item_thm_main_assumption_on_r_k_d}--\ref{item_thm_main_assumption_on_the_other_pair_doubly_stndrd} also
hold for the same $L_1$, $L_2$, $d$, $k$ if we modify $r$ to a smaller number which concludes the proof
also for any $r\leqslant r'$.
\end{prf}

\section{Future projects}
\label{sec_future}
\subsection{Non-divisible line bundle \texorpdfstring{$A$}{A}}
In Theorem~\ref{thm_main_with_splitting} we use a specific splitting $L=A\otimes B$ with $A$ which is
a sufficiently large power of some other line bundle.
One may wonder if it is necessary, if we can just take a splitting $L =A \otimes B$ for sufficiently
ample line bundles $A$ and $B$.
It is very likely that a generalisation of the arguments in this article can provide such result, by considering
triply standard graded $L_1,L_2,L_3$ and appropriately expressing
$A = L_1^{k} \otimes L_2$ and $B=L_1^{l}\otimes L_3$.

\subsection{Apolarity for divisorial algebras and Mori Dream Spaces}

Apolarity theory outlined in \S\ref{sec_apolarity}
for double section rings, is clearly a very special case of more general theory that should be valid for
any divisorial algebras. This generalisation and most importantly its
applications must be carefully and systematically introduced in a follow up work.
In particular, the special case of a projective normal varieties that are Mori Dream Spaces
(their class group and Cox ring are finitely generated) includes homogeneous spaces other than products
of projective spaces.
One should also  compare it to other approaches to apolarity, including
\cite{arrondo_bernardi_marques_mourrain_skew_symmetric_tensor_decomposition}, \cite{staffolani_Schur_apolarity}.

\subsection{Scheme or ideal theoretic equations}

One may ask if the minors obtained in Theorem~\ref{thm_main_intro_simpler}
actually generate the ideal of the cactus variety by analogy to
\cite[Thm~1.1]{sidman_smith_linear_det_eq_for_all_proj_schemes}, or at least if the equations are scheme theoretic.
At the moment we expect this claim to fail. Instead, it is necessary to introduce a potentially non-reduced
scheme structure on $\cactus{r}{X}$ and perhaps the equations are scheme or ideal theoretic for
this enhanced scheme structure.
The article \cite{jabu_keneshlou_cactus_scheme} is the first step in this direction.

\subsection{Explicit conditions on sufficiently ample line bundles}

In specific cases of interest, it could be possible to go through the lines of the proofs in this article
and determine precisely or bound the sufficiently ample condition and describe the set of line bundles $L$
satisfying the conclusions of Theorem~\ref{thm_main_intro_simpler}.
We comment this issue also in \S\ref{sec_ample_and_nef_cones}.

\subsection{Other properties of sufficiently ample line bundles}

Other than the properties mentioned in this article, what are the interesting classes of properties
that are in general satisfied by sufficiently ample line bundles?
In specific cases, such as smooth curves, this topic is being explored. Can we say something more general
than in dimension $1$?

\appendix
\section{Appendix: generalisations to non-algebraically closed fields}
\label{sec_appendix_nonclosed_fields}

So far, throughout the paper, we assumed that the base field $\kk$ is algebraically closed.
In this Appendix we briefly comment on analogous result to Theorem~\ref{thm_main_intro_simpler}
when the base field fails to be algebraically closed.

So assume $\kk$ is an arbitrary field of any characteristics, let $\bar{\kk}$ be its algebraic closure.
Let $X$ be a projective scheme over $\kk$. Denote by $\overline{X} = X\times_{\kk} \bar{\kk}$
the extension of $X$ to a projective scheme over $\bar{\kk}$.
In this setting, the definition of secant or cactus varieties is more tricky than in \S\ref{sec_intro_secant_and_cactus}
and it is discussed in detail in \cite[\S6]{jabu_jelisiejew_finite_schemes_and_secants}.

Explicitly, if $X\subset \PP_{\kk}^N$ then $\sigma_r(X)$ is the smallest subscheme of $\PP_{\kk}^N$
whose extension $\overline{\sigma_r(X)}= \sigma_r(X)\times_{\kk} \bar{\kk}$ contains $\sigma_r(\overline{X})$.
Similarly, $\cactus{r}{X}$ is the smallest subscheme of $\PP_{\kk}^N$  whose
extension $\overline{\cactus{r}{X}}= \cactus{r}{X}\times_{\kk} \bar{\kk}$ contains $\cactus{r}{\overline{X}}$.
See \cite[Sect.~5.6, 5.7, 6.1]{jabu_jelisiejew_finite_schemes_and_secants}
for detailed construction involving Hilbert schemes and relative linear spans of families of schemes
and \cite[Prop.~6.11]{jabu_jelisiejew_finite_schemes_and_secants} for the base change claim.

\begin{thm}
Suppose $X$ is a projective scheme over $\kk$ and $r'$ is a positive integer.
Then there exists a line bundle $\Lambda_0$ on $X$ such that for all ample line bundles $\Lambda$ on $X$
the embedding of $X$ into $\PP\left(H^0(X, \Lambda_0\otimes \Lambda)\right)$
has its cactus variety $\cactus{r}{X}$ set-theoretically defined by $(r+1)\times (r+1)$-minors of a matrix
of linear forms arising from some splitting $L=A \otimes B$ for some line bundles $A$ and $B$ on $X$.
\end{thm}
\begin{prf}[ (sketch)]
Let $\xi\colon \overline{X} \to X$ be the natural morphism coming from $\Spec \bar{\kk} \to \Spec \kk$.
We have an inclusion $\xi^*\colon \Pic (X) \hookrightarrow \Pic(\overline{X})$ by
\cite[Prop.~2.2(i)]{grothendieck_technique_descente_theoremes_d_existence_GA_V_schemas_de_Pic_theoremes_d_exist} and $L\in \Pic(X)$ is ample
if and only if $\xi^*L$ is ample: this equivalence follows from \cite[Thm~4.2(a)]{milne_Abelian_varieties} and
the fact that ampleness can be expressed as an appropriate vanishing of higher cohomologies
of suitable sheaves  \cite[\href{https://stacks.math.columbia.edu/tag/0B5U}{Tag 0B5U}]{stacks_project}.

By Theorem~\ref{thm_main_combinatorial_part}
there exist line bundles $\overline{L}_1$ and $\overline{L}_2$ on $\overline{X}$
such that assumptions of Theorem~\ref{thm_main_with_splitting} are satisfied.
Moreover, by Remark~\ref{rem_sublattice_is_good_enough} we can make sure that these bundles are in
$\xi^*(\Pic (X))$. By Theorem~\ref{thm_main_with_splitting},
the multiplication of sections of $\overline{A}$ and $\overline{B}$
produces a matrix $M$
whose $(r+1)\times (r+1)$ minors define the cactus variety
$\cactus{r}{\overline{X}}$.

By flat base change of cohomology
\cite[\href{https://stacks.math.columbia.edu/tag/02KH}{Lemma 02KH}]{stacks_project} we have
$H^0(X,\overline{A})= H^0(X,A) \otimes_{\kk}\bar{\kk}$ and analogously for $B$ and $L$.
Thus the coefficients of the entries of $M$ are from $\kk$, and so are the coefficients of the minors of $M$.
Therefore, the minors define a subscheme in $\PP(H^0(X, L))$ which after base change
to $\bar{\kk}$ becomes a scheme supported on $\cactus{r}{\overline{X}}$.
By the definition discussed above,
this means that the reduced subscheme of the scheme defined by minors is precisely $\cactus{r}{X}$, concluding the proof.
\end{prf}

\bibliography{cactus-equations.bbl}
\bibliographystyle{alpha}
\end{document}